\documentclass{smfart}
\usepackage[all]{xy}
\usepackage{amssymb}

\newtheorem{theoreme}{Th{\'e}or{\`e}me}[section]
\newtheorem{lemme}{Lemme}[section]
\newtheorem{proposition}{Proposition}[section]
\newtheorem{corollary}{Corollaire}[section]

\theoremstyle{definition}

\theoremstyle{remark}
\newtheorem{remarque}{Remarque}[section]
\numberwithin{equation}{section}

\newcommand{\internalcomment}[1]{}

\newcommand{\kar}{\rule[-0.1mm]{0.12mm}{1.5mm}\hspace{-0.36mm} \ni}
\newcommand{\dgcat}{\mathsf{dgcat}}

\renewcommand{\L}{\mathbf{L}}

\renewcommand{\phi}{\varphi}

\begin{document}
\title{Invariants additifs de dg-cat{\'e}gories}

\author{Gon{\c c}alo Tabuada}
\address{Universit{\'e} Paris 7 - Denis Diderot, UMR 7586 du CNRS, case
  7012, 2, place Jussieu, 75251 Paris cedex 05, France}

\thanks{Soutenu par FCT-Portugal, bourse SFRH/BD/14035/2003}

\date{11 Juillet 2005}

\email{
\begin{minipage}[t]{5cm}
tabuada@math.jussieu.fr
\end{minipage}
}

\begin{abstract}
A l'aide des outils de l'alg{\`e}bre homotopique de Quillen, on  construit `l'invariant additif universel', c'est-{\`a}-dire un
foncteur d{\'e}fini sur la cat{\'e}gorie des
petites dg-cat{\'e}gories et {\`a} valeurs dans une cat{\'e}gorie additive qui rend
inversibles les dg-foncteurs de Morita, transforme les d{\'e}compositions
semi-orthogonales au sens de Bondal-Orlov \cite{Bon-Orl} en sommes directes et qui est
universel pour ces propri{\'e}t{\'e}es. Nous comparons notre construction {\`a}
celle de Bondal-Larsen-Lunts \cite{Grothendieck}.
\end{abstract}
\maketitle

\begin{altabstract}
With the help of the tools of Quillen's homotopical algebra, we construct `the universal additive invariant', namely a
functor from the
category of small dg categories to an additive category, that inverts
the Morita dg functors, transforms the semi-orthogonal
decompositions in the sense of Bondal-Orlov \cite{Bon-Orl} into direct sums and which is
universal for these properties. We compare our construction with that
of Bondal-Larsen-Lunts \cite{Grothendieck}.
\end{altabstract}
\maketitle


\tableofcontents

\section*{Introduction in English} 
In this article, we continue the study \cite{cras} \cite{Toen} of the
category of small differential graded categories ($=$ dg categories)
from the point of view of Quillen model structures
\cite{Quillen}. Concretely, we construct a cofibrantly generated Quillen
model structure on the category of small dg categories such that the
weak equivalences are the {\em Morita} dg functors, i.e. the dg functors $F:\mathcal{A}\rightarrow \mathcal{B}$ that induce an
equivalence $\mathcal{D}(\mathcal{B}) \stackrel{\sim}{\rightarrow}
\mathcal{D}(\mathcal{A})$ between the derived categories. In the
resulting homotopy category $\mathsf{Hmo}$, the derived equivalences
in the sense of \cite{Rickard} \cite{Rickard1} \cite{Keller} correspond
to isomorphisms and the `derived Picard group' of \cite{Rouquier}
can be interpreted as a group of automorphisms.
\par
The category $\mathsf{Hmo}$ is strongly related to the category of
small triangulated categories but carries a richer
structure. It gives us the right framework to formulate universal
properties such as the dg quotient of \cite{Drinfeld}, the
pre-triangulated hull of \cite{Bondal} or the orbit categories of
\cite{Orbit}. 
\par
Another motivation for its study comes from non commutative algebraic
geometry in the sense of Drinfeld \cite{Drinfeldc} and Kontsevich
\cite{Kontsevichc} \cite{Kontsevich}, i.e. the study of dg
categories and their homological invariants. In this spirit, we construct `the universal additive invariant', by
which we mean a functor
$$
\mathcal{U}_{a}:\dgcat \rightarrow \mathsf{Hmo}_0
$$
to an additive category that inverts the Morita dg functors,
transforms the semi-orthogonal decompositions \cite{Bon-Orl} into
direct sums and which is universal for these properties. For instance, $K$-theory and  cyclic homology are
additive invariants and so they factor through $\mathcal{U}_{a}$. We
remark that in $\mathsf{Hmo}_0$, the functor $K_0$ becomes
corepresentable, which immediately yields the Chern characters.
\par
The category $\mathsf{Hmo}_0$ is strongly related to the
Grothendieck ring ${\mathcal{P}\mathcal{T}}_{kar}$ of Karoubian pre-triangulated dg
categories of \cite{Grothendieck}. We make this connexion more precise
by pointing out a
surjection
$$
  {\mathcal{P}\mathcal{T}}^{cl}_{kar}  \rightarrow \mbox{K}_0(\mathsf{Hmo}_0^{cl})
$$ 
after having imposed suitable finiteness conditions.

\section{Introduction}
Dans cet article, nous poursuivons l'{\'e}tude \cite{cras} \cite{Toen} de
la cat{\'e}gorie des petites cat{\'e}gories diff{\'e}rentielles gradu{\'e}es
($=$ dg-cat{\'e}gories) du point de vue des cat{\'e}gories de mod{\`e}les de Quillen
\cite{Quillen}. Plus pr{\'e}cis{\'e}ment, nous munissons la cat{\'e}gorie des
petites dg-cat{\'e}gories d'une structure de cat{\'e}gorie de mod{\`e}les {\`a} engendrement
cofibrant dont les {\'e}quivalences
faibles sont exactement les dg-foncteurs {\em de Morita}, c'est-{\`a}-dire
les dg-foncteurs $F:\mathcal{A} \rightarrow \mathcal{B}$ qui
induisent une {\'e}quivalence $\mathcal{D}(\mathcal{B})
\stackrel{\sim}{\rightarrow} \mathcal{D}(\mathcal{A})$ entre
cat{\'e}gories d{\'e}riv{\'e}es. Dans la cat{\'e}gorie homotopique $\mathsf{Hmo}$
obtenue ainsi, les {\'e}quivalences d{\'e}riv{\'e}es au sens de \cite{Rickard}
\cite{Rickard1} \cite{Keller} correspondent {\`a} des isomorphismes et le `groupe de Picard d{\'e}riv{\'e}' de \cite{Rouquier} y appara{\^\i}t comme un groupe
d'automorphismes. 
\par
La cat{\'e}gorie $\mathsf{Hmo}$ est fortement reli{\'e}e {\`a} la cat{\'e}gorie des
petites cat{\'e}gories triangul{\'e}es mais s'en distingue par sa structure
plus riche. Elle donne un cadre commode pour formuler des propri{\'e}t{\'e}s
universelles comme celles du dg-quotient de \cite{Drinfeld}, de
l'envelope pr{\'e}triangul{\'e}e de \cite{Bondal} ou des cat{\'e}gories d'orbites
de \cite{Orbit}. 
\par
Une autre motivation pour son {\'e}tude provient de la
g{\'e}om{\'e}trie alg{\'e}brique non commutative au sens de Drinfeld
\cite{Drinfeldc} et Kontsevich \cite{Kontsevichc} \cite{Kontsevich},
c'est-{\`a}-dire l'{\'e}tude des dg-cat{\'e}gories et de leurs invariants
homologiques. Dans cette veine, nous construisons `l'invariant additif
universel', c'est-{\`a}-dire un foncteur
$$
\mathcal{U}_a:\dgcat \rightarrow \mathsf{Hmo}_0
$$
{\`a} valeurs dans une cat{\'e}gorie additive qui rend inversibles les
dg-foncteurs de Morita, transforme les d{\'e}compositions
semi-orthogonales \cite{Bon-Orl} en sommes directes et qui est universel pour ces propri{\'e}t{\'e}s. Par exemple, la $K$-th{\'e}orie et
l'homologie cyclique sont des invariants additifs et se factorisent
donc par $\mathcal{U}_a$. Nous observons que dans $\mathsf{Hmo}_0$, le
foncteur $K_0$ devient corepr{\'e}sentable, ce qui donne imm{\'e}diatement les
caract{\`e}res de Chern.
\par
 La cat{\'e}gorie $\mathsf{Hmo}_0$ est {\'e}troitement
reli{\'e}e au anneau de Grothendieck ${\mathcal{P}\mathcal{T}}_{kar}$ des
  dg-cat{\'e}gories pr{\'e}triangul{\'e}es Karobiennes de \cite{Grothendieck}. Nous pr{\'e}cisons
  ce lien en exhibant une surjection
$$
    {\mathcal{P}\mathcal{T}}^{cl}_{kar} \rightarrow  \mbox{K}_0(\mathsf{Hmo}_0^{cl})
$$ 
apr{\`e}s avoir impos{\'e} des conditions de finitude convenables.

\section{Conventions}

Dans toute la suite, $k$ d{\'e}signe un anneau commutatif avec $1$.
Le produit tensoriel $\otimes$ d{\'e}signe toujours le produit tensoriel
sur $k$. Par une {\em dg-cat{\'e}gorie}, nous entendons une $k$-cat{\'e}gorie
diff{\'e}rentielle gradu{\'e}e, voir \cite{Keller} \cite{Drinfeld}.
Soit $\dgcat$ la cat{\'e}gorie des petites dg-cat{\'e}gories. 
Pour la construction du foncteur $\mbox{pre-tr}$, voir \cite{Bondal}.
Pour une dg-cat{\'e}gorie $\mathcal{A}$, on note $\hat{?}:\mathcal{A}
\rightarrow \mathrm{Mod}\,\mathcal{A}$ le dg-foncteur de Yoneda et
$\mathcal{D}(\mathcal{A})$ la cat{\'e}gorie d{\'e}riv{\'e}e, voir \cite{Keller} \cite{Drinfeld}.
Pour les cat{\'e}gories de mod{\`e}les de Quillen, nous renvoyons {\`a}
\cite{Hovey} \cite{Quillen}.

\section{dg-foncteurs quasi-{\'e}quiconiques}

Un dg-foncteur $F$ de $\mathcal{C}$ vers $\mathcal{D}$ est 
{\em quasi-{\'e}quiconique} si:
\begin{itemize}
\item[-] pour tous objets $c_1$ et $c_2$ dans $\mathcal{C}$, le morphisme
  de complexes de $\mathrm{Hom}_{\mathcal{C}}(c_1, c_2)$ vers
  $\mathrm{Hom}_{\mathcal{D}}(F(c_1),F(c_2))$ est un
  quasi-isomorphisme et
\item[-] le foncteur $\mathrm{H}^0(\mbox{pre-tr}(F))$ de
  $\mathrm{H}^0(\mbox{pre-tr}(\mathcal{C}))$ vers
  $\mathrm{H}^0(\mbox{pre-tr}(\mathcal{D}))$ est essentiellement surjectif. 
\end{itemize}
On introduira une structure de cat{\'e}gorie de mod{\`e}les de
Quillen {\`a} engendrement cofibrant dans $\dgcat$ dont les
{\'e}quivalences faibles sont les dg-foncteurs 
quasi-{\'e}quiconiques. Pour cela, on se servira du th{\'e}or{\`e}me~2.1.19
de \cite{Hovey}.
Nous d{\'e}finissons $\mathcal{K}(n)$, $ n \in \mathbb{Z}$, comme
la dg-cat{\'e}gorie avec deux objets $1$, $2$ et dont
les morphismes sont engendr{\'e}s par $f \in \mathrm{Hom}_{\mathcal{K}(n)}^n (1,2)$,
$g \in \mathrm{Hom}_{\mathcal{K}(n)}^{-n} (2,1)$,
$r_1 \in\mathrm{Hom}_{\mathcal{K}(n)}^{-1} (1,1)$,
$r_2 \in \mathrm{Hom}_{\mathcal{K}(n)}^{-1} (2,2)$ et $r_{12}
\in \mathrm{Hom}_{\mathcal{K}(n)}^{n-1} (1,2)$ soumis aux relations $df=dg=0$,
$dr_1=gf-\mathbf{1}_1$, $dr_2 =fg-\mathbf{1}_1$ et $dr_{12}=fr_1 - r_2f$.
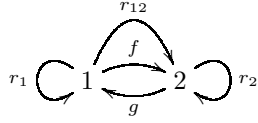
$$\xymatrix{
    1 \ar@(ul,dl)[]_{r_1} \ar@/^/[r]^f \ar@/^0.8cm/[r]^{r_{12}} &
    2 \ar@(ur,dr)[]^{r_2} \ar@/^/[l]^g }
$$
On pose $\mathcal{K}=\mathcal{K}(0)$. Soit $\mathcal{A}$ la dg-cat{\'e}gorie avec un seul object
$3$ et telle que $\mathrm{Hom}_{\mathcal{A}}(3,3)=k$. Soit $F(n)$, $n
\in \mathbb{Z}$, le dg-foncteur de $\mathcal{A}$ vers $\mathcal{K}(n)$ qui envoie $3$
sur $1$.
Pour un entier $n >0$ et des entiers $k_0, \ldots , k_n$, soit
$\mathcal{M}_n(k_0, \ldots , k_n)$ la dg-cat{\'e}gorie avec $n+1$ objets
$0, \ldots , n$ et dont les morphismes sont engendr{\'e}s par des
morphismes $q_{i,j}$ d{\'e}finis pour $0 \leq j < i \leq n$, de source $j$ et de
but $i$, de degr{\'e} $k_i-k_j +1$ et soumis aux relations 
$$d(Q)+Q^2=0\,,$$
o{\`u} $Q$ est la matrice strictement triangulaire inf{\'e}rieure de
coefficients $q_{i,j}$ pour $0 \leq j < i \leq n$.

Soit $\mathrm{cone}_n(k_0, \ldots , k_n)$ la sous-dg-cat{\'e}gorie pleine de la
dg-cat{\'e}gorie des dg-modules ({\`a} droite) sur $\mathcal{M}_n(k_0, \ldots , k_n)$
dont les objets sont les $\hat{l}$, $0 \leq l \leq n$, et le {\em c{\^o}ne
it{\'e}r{\'e}} $X_n$, qui a le m{\^e}me module gradu{\'e} sous-jacent que le dg-module
$$ X=\bigoplus^n_{l=0} \hat{l}\left[ k_l \right]$$
et dont la diff{\'e}rentielle est $d_X+\widehat{Q}$. 

Soit $L:\mathcal{A} \rightarrow \mathrm{cone}_n(k_0, \ldots , k_n)$ le dg
foncteur qui envoie $3$ sur $X_n$. On consid{\`e}re la somme amalgam{\'e}e
$$\xymatrix{
\mathcal{A} \ar[r]^-{L} \ar[d]_F \ar@{}[dr]|{\lrcorner} & \mathrm{cone}_n(k_0, \ldots , k_n) \ar[d] \\
\mathcal{K} \ar[r] & \mathrm{cone}_n(k_0, \ldots ,
k_n) \amalg_{\mathcal{A}} \mathcal{K} 
}$$
et on d{\'e}finit $\mathrm{coneh}_n(k_0, \ldots , k_n)$ comme la sous-dg-cat{\'e}gorie
pleine de la somme amalgam{\'e}e dont les objets sont les images des
objets $\hat{l}$, $0 \leq l \leq n$, de $\mathrm{cone}_n(k_0, \ldots , k_n)$
et l'image $X\!h\,_n$ de l'object $2$ de $\mathcal{K}$. On note $I_n(k_0, \ldots , k_n)$ le dg foncteur fid{\`e}le, mais non plein, de 
$\mathcal{M}_n(k_0, \ldots , k_n)$ dans $\mathrm{coneh}_n(k_0, \ldots , k_n)$.

\begin{theoreme}\label{theorem}
 Si on consid{\`e}re pour cat{\'e}gorie $\mathcal{C}$ la
  cat{\'e}gorie $\dgcat$, pour classe $W$ la sous-cat{\'e}gorie de $\dgcat$ des
dg-foncteurs quasi-{\'e}quiconiques, pour classe $J$ les dg-foncteurs $F$ et $R(n)$,
$n\in \mathbb{Z}$, (voir \cite{cras}), $F(n)$, $n \in \mathbb{Z}$, et $I_n(k_0,
 \ldots , k_n)$ et pour classe $I$ les dg-foncteurs $Q$ et $S(n)$,
$n\in \mathbb{Z}$, (voir \cite{cras}) alors les conditions du th{\'e}or{\`e}me \cite[2.1.19]{Hovey} sont satisfaites.
\end{theoreme}

\begin{remarque}\label{fibrants}
On peut montrer ais{\'e}ment que pour la structure obtenue, les objets
fibrants sont les dg-cat{\'e}gories $\mathcal{A}$, telles que l'image
essentielle du plongement $\mathrm{H}^0(\mathcal{A}) \hookrightarrow
\mathcal{D}(\mathcal{A})$ est stable par suspension et c{\^o}nes. 
Cela est {\'e}quivalent {\`a} ce que le dg-foncteur $\mathcal{A} \hookrightarrow
\mbox{pre-tr}(\mathcal{A})$, voir \cite{Bondal}, soit une quasi-{\'e}quivalence. 
\end{remarque}

On observe facilement que les conditions $\mbox{{\it (i)}}$,
$\mbox{{\it (ii)}}$ et $\mbox{{\it (iii)}}$ du th{\'e}or{\`e}me~2.1.19
de \cite{Hovey} sont verifi{\'e}es.

D{\'e}montrons maintenant que   $J-\mbox{inj}\cap W = I-\mbox{inj}$.
Pour cela, on consid{\`e}re la classe $\mbox{{\bf Surj}}$ form{\'e}e des foncteurs
$G : \mathcal{H} \rightarrow \mathcal{I}$ dans $\dgcat$ qui
v{\'e}rifient~:
\begin{itemize}
\item[-] $G$ induit une surjection de l'ensemble des objets de $\mathcal{H}$ sur l'ensemble des
objets de $\mathcal{I}$ et
\item[-] $G$ induit des quasi-isomorphismes surjectifs
dans les complexes de morphismes.
\end{itemize}

\begin{lemme}\label{I-inj-Surj}
On a $\mbox{{\bf Surj}} = I-\mbox{inj}$ (voir \cite {cras}).
\end{lemme}

\begin{lemme}\label{J-inj-Surj 1}
On a $\mbox{{\bf Surj}} \supseteq J-\mbox{inj} \cap W$.
\end{lemme}
\begin{proof}
Soit $L$ un dg-foncteur de $\mathcal{D}$ vers
$\mathcal{S}$ qui appartient {\`a} $J-\mbox{inj}\cap W$. La classe
$R(n)-\mbox{drt}$ est form{\'e}e des dg-foncteurs surjectifs au niveau des
complexes de morphismes. Comme $L \in W$, il suffit de montrer que $L$
est surjectif au niveau des objets. Soit $E \in \mathcal{S}$ un
object quelconque. Comme $L \in W$, il existe $C \in
\mbox{pre-tr}(\mathcal{D})$ et un morphisme ferm{\'e} $q$ de
$\mbox{pre-tr}(\mathcal{S})$ qui fournit un isomorphisme entre l'image
par $\mbox{pre-tr}(L)$ de $C$ et $E$ dans
$\mathrm{H}^0(\mbox{pre-tr}(\mathcal{S}))$.  On identifie les
dg-cat{\'e}gories $\mathcal{D}$ et $\mathbb{Z}\,\mathcal{D}$ avec leurs
images dans $\mbox{pre-tr}(\mathcal{D})$. Il existe alors trois possibilit{\'e}s~:
\begin{enumerate}
\item l'objet $C$ appartient {\`a} la dg-cat{\'e}gorie $\mathcal{D}$. Alors on est dans les conditions de
  \cite{cras}.
\item l'objet $C$ appartient {\`a} la dg-cat{\'e}gorie $\mathbb{Z}\,\mathcal{D}$. Alors $C$ est de la forme
  $C=(\overline{C}, n)$, o{\`u} $ \overline{C} \in \mathcal{D}$ et $n \in \mathbb{Z}$.
On a la situation suivante
$$\xymatrix{
C=((\overline{C}), n) \ar@{|->}[d]^{\mathbb{Z}(L)} & \\
\mathbb{Z}(L)(C)=(L(\overline{C}), n) \ar[r]^-{q} & (E,0).
}$$
Ainsi, $q$ est l'image de $f$ par un dg-foncteur de $\mathcal{K}(-n)$
vers $\mathcal{S}$. Comme on a $L \in J-\mbox{inj}$, on peut relever
le morphisme $q$ et par cons{\'e}quence l'objet $E$.
\item l'objet $C$ appartient {\`a} la dg-cat{\'e}gorie
  $\mbox{pre-tr}(\mathcal{D})$ mais non pas {\`a} la dg-cat{\'e}gorie
  $\mathbb{Z}\,\mathcal{D}$. On sait d'apr{\`e}s \cite{Bondal} que, dans
  $\mbox{pre-tr}(\mathcal{D})$, l'objet $C$ s'{\'e}crit d'une fa{\c c}on
  canonique comme un c{\^o}ne it{\'e}r{\'e} sur des morphismes de
  $\mathcal{D}$. Comme le dg-foncteur $\mbox{pre-tr}(L)$ pr{\'e}serve les c{\^o}nes, l'object
  $\mbox{pre-tr}(L)(C)$ s'identifie au c{\^o}ne it{\'e}r{\'e} sur leurs images par
  $L$. On peut donc construire le carr{\'e} commutatif suivant
$$\xymatrix @R=1pc{
*+<1pc>{\mathcal{M}_n(k_0, \ldots , k_n)} \ar[r] \ar@{^{(}->}[dd]_{can}
& \mathcal{D} \ar[d]^L\\
 & \mathcal{S} \ar[d] \\
\mathrm{cone}_n(k_0, \ldots , k_n) \ar[r]^-H & \mbox{pre-tr}(\mathcal{S})\,.
}$$
Ainsi, $q$ est l'image de $f$ par un dg-foncteur de $\mathcal{K}$ vers
$\mbox{pre-tr}(\mathcal{S})$ qui envoie l'object $1$ sur $H(X_n)$ et l'object $2$
sur l'object $E$. Le dg-foncteur $H$ s'{\'e}tend donc en un dg-foncteur
$\overline{H}$ de $\mathrm{cone}_n(k_0, \ldots , k_n) \amalg_{\mathcal{A}}
\mathcal{K}$ vers $\mbox{pre-tr}(\mathcal{S})$. On observe que la
restriction de $\overline{H}$, qu'on note $\widetilde{H}$, {\`a} la dg-cat{\'e}gorie
pleine $\mathrm{coneh}_n(k_0, \ldots , k_n)$, a son image dans la dg-cat{\'e}gorie
$\mathcal{S}$. Cela permet de construire le diagramme commutatif suivant 
$$\xymatrix@!0 @R=3pc @C=4pc{
 & & \mathcal{M}_n \ar[rr] \ar[dr]_-{I_n} \ar[dl]^-{can} & 
 & \mathcal{D} \ar[d]^L \\
\mathcal{A} \ar@{}[drr]|{\lrcorner} \ar[r] \ar[dr]_-F & \mathrm{cone}_n \ar[dr] \ar[drrr]^(0.3){H} & &
 \mathrm{coneh}_n \ar[r]^-{\widetilde{H}}
 \ar@{^{(}->}[dl]|\hole \ar@{.>}[ur]  & \mathcal{S} \ar[d] \\
 & \mathcal{K} \ar[r] & \mathrm{cone}_n \amalg_{\mathcal{A}} \mathcal{K} \ar[rr]_-{\overline{H}} &  & \mbox{pre-tr}(\mathcal{S})\,.
}$$
Comme $L$ appartient {\`a} $J-\mbox{inj}$, l'objet $E$ est bien l'image d'un objet de $\mathcal{D}$.
\end{enumerate}
\end{proof}

Soit $\mathcal{B}$ une dg-cat{\'e}gorie et $C$ un dg-foncteur de $\mathcal{M}_n(k_0, \ldots , k_n)$ vers $\mathcal{B}$.
\begin{remarque}
Le dg-foncteur $C$ se prolonge en un dg-foncteur
$$ C^{\ast}: \mathrm{Mod}\,\mathcal{M}_n(k_0, \ldots , k_n)
\rightarrow \mathrm{Mod}\,\mathcal{B}$$
adjoint {\`a} gauche de la restriction le long de $C$. Le dg-foncteur $C$
{\em admet un c{\^o}ne} $Y$ dans $\mathcal{B}$ si $C^{\ast}(X_n)$ est
repr{\'e}sentable par un objet $Y$ de $\mathcal{B}$. C'est le cas ssi le
foncteur $C$ se prolonge en un dg-foncteur
$$\widetilde{C}: \mathrm{cone}_n(k_0, \ldots , k_n) \rightarrow \mathcal{B}\,.$$
Si $G: \mathcal{B} \rightarrow \mathcal{B}'$ est un dg-foncteur et
$C$ admet le c{\^o}ne $Y$ dans $\mathcal{B}$, alors $G\circ C$ admet le
c{\^o}ne $G(Y)$ dans $\mathcal{B}'$.
Si $C$ se factorise par la sous-dg-cat{\'e}gorie pleine
$\mathcal{B}\,\backslash \{Y \}$ et admet le c{\^o}ne $Y$ dans
$\mathcal{B}$, alors le carr{\'e}
$$\xymatrix{
*+<1pc>{\mathcal{M}_n(k_0, \ldots , k_n)} \ar[r] \ar@{^{(}->}[d] \ar@{}[dr]|{\lrcorner} & *+<1pc>{\mathcal{B}\,\backslash
\{Y \}} \ar@{^{(}->}[d]\\
\mathrm{cone}_n(k_0, \ldots , k_n) \ar[r] & \mathcal{B}
}$$
est cocart{\'e}sien.
\end{remarque}

Ces remarques impliquent le lemme suivant.
\begin{lemme} \label{lemme-clef}
Le carr{\'e} commutatif
$$\xymatrix{
*+<1pc>{\mathcal{M}_n(k_0, \ldots , k_n)} \ar[rr]^-{I_n(k_0, \ldots , k_n)}
\ar@{^{(}->}[d]_-{can}  \ar@{}[drr]|{\lrcorner} & &  *+<1pc>{\mathrm{coneh}_n(k_0, \ldots , k_n)}
\ar@{^{(}->}[d] \\
\mathrm{cone}_n(k_0, \ldots , k_n) \ar[rr]  & &  \mathrm{cone}_n(k_0, \ldots , k_n) \amalg_{\mathcal{A}}\mathcal{K}  
}$$
est cocart{\'e}sien.
\end{lemme}

\begin{lemme}\label{J-inj-Surj 2}
On a $\mbox{{\bf Surj}} \subseteq J-\mbox{inj} \cap W$.
\end{lemme}
\begin{proof}
Soit $H$ un dg-foncteur
de $\mathcal{N}$ vers $\mathcal{E}$ dans la classe $\mbox{{\bf Surj}}$. Comme $H$ est surjectif au niveau des objets et un
quasi-isomorphisme au niveau des complexes de morphismes, on a $H \in
W$. La classe $R(n)-\mbox{drt}$ est form{\'e}e des dg-foncteurs surjectifs
au niveau des complexes de morphismes. Il suffit de montrer que $H$ a
la propri{\'e}t{\'e} de rel{\`e}vement {\`a} droite par rapport {\`a} $F$, $F(n)$, $n \in
\mathbb{Z}$, et $I_n(k_0, \ldots , k_n)$. On consid{\`e}re ces trois cas~:
\begin{enumerate}
\item on sait par \cite{cras} que $H \in F-\mbox{drt}$;
\item en effet, on a $H \in F(n)-\mbox{drt}$, $n \in \mathbb{Z}$, par un argument compl{\`e}tement
  analogue au cas pr{\'e}c{\'e}dent;
\item on consid{\`e}re le diagramme suivant~:
$$\xymatrix{
\mathcal{M}_n(k_0,\ldots, k_n) \ar[r] \ar[d]_{I_n(k_0,
  \ldots , k_n)} & \mathcal{N} \ar[d]^H\\
\mathrm{coneh}_n(k_0, \ldots , k_n) \ar[r]_-{L} & \mathcal{E}
}$$
Le lemme~\ref{lemme-clef} permet d'{\'e}tendre le dg-foncteur $L$ en un dg-foncteur
$\overline{L}$ de $\mathrm{cone}_n(k_0, \ldots , k_n)\amalg_{\mathcal{A}}
 \mathcal{K}$ vers $\mbox{pre-tr}(\mathcal{E})$. 
$$\xymatrix@!0 @R=2pc @C=4pc{
&  & \mathcal{M}_n \ar[rr] \ar[dl] \ar[dd]|\hole & &  *+<1pc>{\mathcal{N}}
\ar[dd]^-H \ar@{^{(}->}[dl] \\
\mathcal{A} \ar@{}[ddr]|{\lrcorner} \ar[dd]_F \ar[r] & \mathrm{cone}_n \ar[dd] \ar[rr]  & &
\mbox{pre-tr}(\mathcal{N}) \ar[dd] & \\
& & \mathrm{coneh}_n \ar[rr]|\hole^(0.3){L} \ar[dl] &    & *+<1pc>{\mathcal{E}}
\ar@{^{(}->}[dl] \\
\mathcal{K} \ar[r] & \mathrm{cone}_n \amalg_{\mathcal{A}} \mathcal{K}
\ar[rr]^-{\overline{L}} &   & \mbox{pre-tr}(\mathcal{E}) & 
}$$
Puisque le dg-foncteur $H$ appartient {\`a} la classe $\mbox{{\bf Surj}}$, on peut
relever le dg-foncteur $\overline{L}$ vers
$\mbox{pre-tr}(\mathcal{N})$. On remarque que la restriction de ce
rel{\`e}vement {\`a} la sous-dg-cat{\'e}gorie pleine $\mathrm{coneh}_n(k_0, \ldots
,k_n)$ fournit un rel{\`e}vement du dg-foncteur $L$ vers la
dg-cat{\'e}gorie $\mathcal{N}$.
\end{enumerate}
\end{proof}

\begin{lemme} \label{J-cell-dans-W} On a $J-\mbox{cell}\subset W$.
\end{lemme}
\begin{proof}
On sait d{\'e}j{\`a} par \cite{cras} que les classes $F-\mbox{cell}$ et
$R(n)-\mbox{cell}$ sont form{\'e}es des quasi-{\'e}quivalences et donc
contenues dans la classe $W$. Soit $n \in \mathbb{Z}$ et $T:\mathcal{A}
\rightarrow \mathcal{J}$ un dg-foncteur quelconque. On consid{\`e}re la
somme amalgam{\'e}e suivante
$$\xymatrix{
\mathcal{A} \ar[d]_{F(n)} \ar[r]^T \ar@{}[dr]|{\lrcorner} & \mathcal{J} \ar[d]^{inc}\\
\mathcal{K}(n) \ar[r] & \mathcal{U} 
}$$ 
dans $\dgcat$. Il s'agit de v{\'e}rifier que $inc \in W$. Il
faut verifier deux conditions:
\begin{enumerate}
\item L'inclusion $\mathrm{Hom}_{\mathcal{J}}(X,Y) \rightarrow \mathrm{Hom}_{\mathcal{U}} (inc(X), inc(Y))$ est un quasi-isomorphisme
  pour tous $X,Y \in \mathcal{J}$. Pour cela, on raisonne comme dans
  le cas o{\`u} on a le dg-foncteur $F$ au lieu de
  $F(n)$. Voir \cite{cras}.
\item Le foncteur
  $\mathrm{H}^0(\mathbb{Z}(inc)):\mathrm{H}^0(\mathbb{Z}\,\mathcal{J})
  \rightarrow \mathrm{H}^0(\mathbb{Z}\,\mathcal{U})$ est essentiellement
  surjectif. En effet, la dg-cat{\'e}gorie $\mathcal{U}$ s'obtient {\`a} partir de
  $\mathcal{J}$ en rajoutant un nouvel objet $S$ homotopiquement {\'e}quivalent {\`a} un d{\'e}cal{\'e} de l'objet $T(3)$. Plus
  pr{\'e}cis{\'e}ment, l'object $(S,0)$
  devient isomorphe {\`a} $(T(3),-n)$ dans $\mathrm{H}^0(\mathbb{Z}\,\mathcal{U})$.
\end{enumerate}

Soit maintenant $T:\mathcal{M}_n(k_0,\ldots, k_n) \rightarrow \mathcal{J}$ un
dg-foncteur quelconque. On consid{\`e}re la somme amalgam{\'e}e suivante
$$\xymatrix{
\mathcal{M}_n(k_0,\ldots, k_n) \ar[d]_{I_n(k_0, \ldots ,
  k_n)} \ar[r]^-{T} \ar@{}[dr]|{\lrcorner}   & \mathcal{J} \ar[d]^{inc} \\
\mathrm{coneh}_n(k_0, \ldots , k_n) \ar[r] & \mathcal{U}
}$$ 
dans $\dgcat$. Il s'agit de v{\'e}rifier que $inc \in W$. Il
faut verifier deux conditions:
\begin{enumerate}
\item L'inclusion $\mathrm{Hom}_{\mathcal{J}}(X,Y) \rightarrow \mathrm{Hom}_{\mathcal{U}} (inc(X), inc(Y))$ est un quasi-isomor\-phisme
  pour tous $X,Y \in \mathcal{J}$. Dans le diagramme commutatif
  suivant, la dg-cat{\'e}gorie $\mathcal{N}$ est la somme amalgam{\'e}e de
  $\mathcal{J}$ et $\mathrm{cone}_n$, et $\mathcal{E}$ est la somme
  amalgam{\'e}e de $\mathcal{N}$ et $\mathcal{K}$.
$$\xymatrix@!0 @R=2pc @C=4pc{
&  & \mathcal{M}_n \ar[rr]^-T \ar[dl] \ar[dd]|\hole & &  \mathcal{J}
\ar[dd]^-{inc} \ar[dl]^L \\
\mathcal{A} \ar@{}[ddr]|{\lrcorner} \ar[dd]_F \ar[r] & \mathrm{cone}_n \ar[dd] \ar[rr] & &
\mathcal{N} \ar[dd]^(0.3){H}  & \\
& & \mathrm{coneh}_n \ar[rr]|\hole \ar@{^{(}->}[dl] &    & *+<1pc>{\mathcal{U}}
\ar@{^{(}->}[dl] \\
\mathcal{K} \ar[r] & \mathrm{cone}_n \amalg_{\mathcal{A}} \mathcal{K}
\ar[rr] &   & \mathcal{E} & 
}$$
Notons que le dg-foncteur $inc$ s'identifie {\`a} la
 composition de $L$ et $H$ suivie d'une restriction. On peut identifier
 la dg-cat{\'e}gorie $\mathcal{N}$ {\`a} la sous-dg-cat{\'e}gorie pleine des dg-modules ({\`a}
 droite) sur $\mathcal{J}$ dont les objets sont les dg-modules
 represent{\'e}es et le $n$-i{\`e}me c{\^o}ne it{\'e}r{\'e} sur l'image de la famille
 $T(q_{i,j})$ par le foncteur de Yoneda. Le dg-foncteur $L$ s'identifie alors au plongement
 de Yoneda et il est donc pleinement fid{\`e}le. Le dg-foncteur $H$ s'identifie au dg-foncteur de $\mathcal{N}$ vers
 la somme amalgam{\'e}e avec $\mathcal{K}$. On sait d'apr{\`e}s \cite{cras}
 que $H$ est donc un quasi-isomorphisme. Cela implique que $inc$
 satisfait la condition.
\item  Le foncteur
  $\mathrm{H}^0(\mbox{pre-tr}(inc)):\mathrm{H}^0(\mbox{pre-tr}(\mathcal{J}))
  \rightarrow \mathrm{H}^0(\mbox{pre-tr}(\mathcal{U}))$ est essentiellement
  surjectif. En effet, la dg-cat{\'e}gorie $\mathcal{U}$ poss{\`e}de un object
  de plus que la dg-cat{\'e}gorie $\mathcal{J}$, {\`a} savoir l'image de $X\!h\,_n$. Soit $T(X_n)$,
  dans $\mbox{pre-tr}(\mathcal{J})$, le $n$-i{\`e}me c{\^o}ne it{\'e}r{\'e} associ{\'e} au
  dg-foncteur $T$. Alors le lemme~\ref{lemme-clef} nous montre que l'image de $T(X_n)$ par le dg-foncteur
  $\mbox{pre-tr}(inc)$ et $X\!h\,_n$ sont homotopiquement {\'e}quivalents
  dans $\mbox{pre-tr}(\mathcal{U})$.
\end{enumerate}
\end{proof} 

Nous avons v{\'e}rifi{\'e} que $J-\mbox{cell}\subseteq W$ (lemme~\ref{J-cell-dans-W}) et
que $I-\mbox{inj}$ est {\'e}gal {\`a} $J-\mbox{inj}\cap W$
(lemmes~\ref{I-inj-Surj}, \ref{J-inj-Surj 1} et \ref{J-inj-Surj2}).
Ces conditions impliquent celles du th{\'e}or{\`e}me de Hovey \cite[2.1.19]{Hovey}.

On note $\mathcal{Q}eq$ la classe des quasi-{\'e}quivalences, voir
\cite{cras}, et $\mathcal{Q}ec$ la classe plus large des
dg-foncteurs quasi-{\'e}quiconiques. On note $\mathsf{Heq}$ et
$\mathsf{Hec}$ les cat{\'e}gories homotopiques de $\dgcat$ par rapport {\`a}
ces classes. Soient $\mathcal{A}$ et $\mathcal{B}$ des dg-cat{\'e}gories et $\mbox{can}_1:\mathsf{Heq} \rightarrow \mathsf{Hec}$ le
foncteur canonique. 

\begin{corollary} \label{adjonct1}
On a une adjonction
$$
\mathrm{Hom}_{\mathsf{Hec}}(\mbox{can}_1(\mathcal{A}),
\mathcal{B}) \stackrel{\sim}{\longrightarrow} \mathrm{Hom}_{\mathsf{Heq}}(\mathcal{A},
\mathcal{B}_{fib})\,,
$$
o{\`u} $\mathcal{B}_{fib}$ est un remplacement fibrant
de la dg-cat{\'e}gorie $\mathcal{B}$ par rapport {\`a} la structure de
cat{\'e}gorie de mod{\`e}les de Quillen quasi-{\'e}quiconique d{\'e}crite
dans th{\'e}or{\`e}me~\ref{theorem}.
\end{corollary}
\begin{proof}
En effet les structures de cat{\'e}gorie de mod{\`e}les de Quillen dans
$\dgcat$ de \cite{cras} et du th{\'e}or{\`e}me~\ref{theorem} ont les m{\^e}mes
cofibrations et les m{\^e}mes fibrations acycliques. Cela implique que les
objets cofibrants et les objets cylindres sont les m{\^e}mes et par
cons{\'e}quent la relation d'homotopie est la m{\^e}me dans les deux situations.
On a donc des bijections
$$
\mathrm{Hom}_{\mathsf{Hec}}(\mbox{can}_1(\mathcal{A}), \mathcal{B})
 \stackrel{\sim}{\longrightarrow}  \mathrm{Hom}_{\dgcat}(\mathcal{A}_{cof},\mathcal{B}_{fib})/ \mbox{htp} \stackrel{\sim}{\longleftarrow} \mathrm{Hom}_{\mathsf{Heq}}(\mathcal{A}, \mathcal{B}_{fib})\,.
$$
\end{proof}

\begin{lemme}\label{homotopes}
Les dg-cat{\'e}gories $\mathcal{B}_{fib}$ et
$\mbox{pre-tr}(\mathcal{B})$ sont isomorphes dans $\mathsf{Heq}$.
\end{lemme}
\begin{proof}
On consid{\`e}re le carr{\'e} commutatif
$$
\xymatrix{
*+<2pc>{\mathcal{B}} \ar[r] \ar@{>->}[d]_-F & \mbox{pre-tr}(\mathcal{B})
\ar[d]^-{pre-tr(F)} \\
*+<1pc>{\mathcal{B}_{fib}} \ar@{^{(}->}[r] &
\mbox{pre-tr}(\mathcal{B}_{fib})\,.
}
$$
Le dg-foncteur $F$ appartient {\`a} la classe $\mathcal{Q}ec$ et donc par
definition $\mbox{pre-tr}(F)$ appartient {\`a} la classe
$\mathcal{Q}eq$. Puisque la dg-cat{\'e}gorie $\mathcal{B}_{fib}$ est
fibrante, la remarque~\ref{fibrants} implique que le dg-foncteur $\mathcal{B}_{fib}
\hookrightarrow \mbox{pre-tr}(\mathcal{B}_{fib})$ appartient {\`a} la
classe $\mathcal{Q}eq$. Cela implique le lemme.
\end{proof}

On note $\mathsf{Heq}_{ex}$ la sous-cat{\'e}gorie pleine de
$\mathsf{Heq}$ dont les objets sont les dg-cat{\'e}gories exactes,
voir \cite[2]{Exact}. On note $\mbox{inc}: \mathsf{Heq}_{ex} \hookrightarrow
\mathrm{Heq}$ l'inclusion. Le corollaire~\ref{adjonct1} et le
lemme~\ref{homotopes} impliquent

\begin{corollary}\label{adjonct3}
On a une adjonction
$$
\mathrm{Hom}_{\mathsf{Heq}}(\mathcal{A}, \mbox{inc}(\mathcal{B})) \stackrel{\sim}{\longleftarrow} 
\mathrm{Hom}_{\mathsf{Heq}}(\mbox{pre-tr}(\mathcal{A}), \mathcal{B})\,.
$$
\end{corollary}
On remarque que les corollaires~\ref{adjonct1} et~\ref{adjonct3} nous
fournissent une {\'e}quivalence de cat{\'e}gories entre $\mathsf{Hec}$ et
$\mathsf{Heq}_{ex}$. On a donc le diagramme suivant
$$
\xymatrix@!0 @R=4pc @C=8pc{
*+<0.5pc>{\mathsf{Heq}} \ar@<-1ex>[d]_-{can_1} \ar@<-1ex>[r]_-{pre-tr} &
*+<1pc>{\mathsf{Heq}_{ex}} 
\ar@<2.5ex>[dl]^-{\sim} \ar@<-1ex>@{^{(}->}[l]_-{inc}\\
*+<0.5pc>{\mathsf{Hec}} \ar@<-1ex>[u]_-{pre-tr} \ar@<-1.5ex>[ur]^-{\sim} & 
}
$$
On note $\mathsf{rep}_{ec}(\mathcal{A}, \mathcal{B})$ la sous-cat{\'e}gorie pleine
de la cat{\'e}gorie d{\'e}riv{\'e}e $\mathcal{D}(\mathcal{A}^{op} \otimes^{\mathbb{L}} \mathcal{B})$,
  voir \cite{Keller} \cite{Drinfeld}, dont les objets sont les dg
  $\mathcal{A}\mbox{-}\mathcal{B}$-bimodules $X$ tels que $X(?, A)$ est
  isomorphe dans $\mathcal{D}(\mathcal{B})$ {\`a} un objet de l'image du
  dg-foncteur canonique $\mbox{pre-tr}(\mathcal{B})\rightarrow
  \mathrm{Mod}\, \mathcal{B}$, pour tout $A \in \mathcal{A}$. On note
  $\left[X \right]$ la classe d'isomorphisme de $X$ dans
  $\mathsf{rep}_{ec}(\mathcal{A},\mathcal{B})$. Pour une
  cat{\'e}gorie essentiellement petite $\mathcal{C}$, on note
  $\mathrm{Iso}(\mathcal{C})$ l'ensemble de ses classes d'isomorphisme.
\begin{corollary}\label{adjonct4}
On a une bijection
$$\mathrm{Hom}_{\mathsf{Hec}}(\mathcal{A},
\mathcal{B}) \stackrel{\sim}{\longrightarrow} \mathrm{Iso}(\mathsf{rep}_{ec}(\mathcal{A}, \mathcal{B}))\,.
$$
\end{corollary}

\begin{proof}
D'apr{\`e}s le corollaire~\ref{adjonct1} et le lemme~\ref{homotopes}, $\mathrm{Hom}_{\mathsf{Hec}}(\mathcal{A},
\mathcal{B})$ s'identifie {\`a}  $\mathrm{Hom}_{\mathsf{Heq}}(\mathcal{A},
\mbox{pre-tr}(\mathcal{B}))$. Par \cite{Toen} on sait que
cet ensemble s'identifie {\`a} $\mbox{Iso}(\mathsf{rep}(\mathcal{A},\mbox{pre-tr}(\mathcal{B}))$.
On remarque que $\mathsf{rep}(\mathcal{A},\mbox{pre-tr}(\mathcal{B}))$
s'identifie {\`a} $\mathsf{rep}_{ec}(\mathcal{A},\mathcal{B})$, d'o{\`u} le resultat.
\end{proof}

\begin{remarque}\label{monoidale}
On sait que les structures de cat{\'e}gorie de mod{\`e}les de Quillen dans
$\dgcat$ de \cite{cras} et du th{\'e}or{\`e}me~\ref{theorem} ont les m{\^e}mes
cofibrations et les m{\^e}mes fibrations acycliques. Cela a pour
consequence que
\begin{itemize}
\item[-] Le foncteur de remplacement cofibrant simplicial $\Gamma^{\ast}$, voir
\cite{Hovey}, dans $\dgcat$ est le m{\^e}me dans les deux situations. Cela
implique que les corollaires~\ref{adjonct1}, ~\ref{adjonct3} et~\ref{adjonct4} sont encore vrais si on remplace $\mathrm{Hom}$ par
l'espace de morphismes, $\underline{\mathrm{Map}}$, voir \cite{Hovey}, et
$\mbox{Iso}(\mathsf{rep}_{ec}(\mathcal{A},\mathcal{B}))$,
dans le corollaire~\ref{adjonct4}, par le nerf de la sous-cat{\'e}gorie de
$\mathrm{Mod}(\mathcal{A}^{op} \otimes^{\mathbb{L}}\mathcal{B})$ dont les objets
sont les m{\^e}mes que ceux de $\mathsf{rep}_{ec}(\mathcal{A},
\mathcal{B})$ et dont les morphismes sont les quasi-isomorphismes.
\item[-] Le produit tensoriel $\otimes^{\mathbb{L}}$ de $\mathsf{Heq}$
  descend {\`a} la localis{\'e}e $\mathsf{Hec}$. La cat{\'e}gorie mono{\"\i}dale sym{\'e}trique $(\mathsf{Hec},\,-\otimes^{\mathbb{L}}-)$ est ferm{\'e}e, voir \cite{Hovey}, et l'espace de
  morphismes interne de $\mathsf{Hec}$, $\mathbb{R}\underline{\mathrm{Hom}}_{\mathsf{Hec}}(\mathcal{A},
    \mathcal{B})$ s'identifie {\`a} $\mathbb{R}\underline{\mathrm{Hom}}_{\mathsf{Heq}}(\mathcal{A},
      \mbox{pre-tr}(\mathcal{B}))$, voir \cite{Toen}.
\end{itemize} 
\end{remarque}

\section{Additivisation}
Soient $\mathcal{A}$, $\mathcal{B}$ et $\mathcal{C}$ des
dg-cat{\'e}gories. La composition dans $\mathsf{Hec}$ est induite par le
bifoncteur
$$
\xymatrix@!0 @R=1,5pc @C=12pc{
-\otimes^{\mathbb{L}}_{\mathcal{B}}-: \mathsf{rep}_{ec}(\mathcal{A},
 \mathcal{B}) \times \mathsf{rep}_{ec}(\mathcal{B}, \mathcal{C})
\ar[r]   & \mathsf{rep}_{ec}(\mathcal{A}, \mathcal{C})\\
*+<2pc>{(X , Y)} \ar@{|->}[r]  & X_{cof}\otimes_{\mathcal{B}}Y\, ,
}
$$
o{\`u} $X_{cof}$ est un remplacement cofibrant de $X$ dans la cat{\'e}gorie
$\mathrm{Mod}(\mathcal{A}^{op}\otimes^{\mathbb{L}}\mathcal{B})$,
par rapport {\`a} sa structure de cat{\'e}gorie de mod{\`e}les de Quillen, voir \cite{Toen}.

\begin{remarque}\label{bitri}
Le bifoncteur $-\otimes^{\mathbb{L}}-$ est bi-triangul{\'e}, puisqu'il est induit par le dg-bifoncteur produit tensoriel de bimodules.
\end{remarque}

Soit $\mathsf{Hec}_0$ la cat{\'e}gorie qui a pour objets les petites
dg-cat{\'e}gories et telle que $\mathrm{Hom}_{\mathsf{Hec}_0}(\mathcal{A},
\mathcal{B})$ est le groupe de Grothendieck de la cat{\'e}gorie triangul{\'e}e
$\mathsf{rep}_{ec}(\mathcal{A}, \mathcal{B})$. L'operation de
composition est induite par le bifoncteur $-\otimes^{\mathbb{L}}-$
d'apr{\`e}s la remarque~\ref{bitri}. On dispose d'un foncteur canonique
$\mbox{add}_1:\mathsf{Hec}\rightarrow \mathsf{Hec}_0$.

\begin{lemme}\label{addit}
La cat{\'e}gorie $\mathsf{Hec}_0$ est additive et le foncteur canonique
$\mathsf{Hec}\rightarrow \mathsf{Hec}_0$ transforme les produits finis
et les coproduits finis en sommes directes.
\end{lemme} 
\begin{proof}
Par construction, les ensembles de morphismes de $\mathsf{Hec}_0$ sont
des groupes ab{\'e}liens et la composition est $\mathbb{Z}$-bilin{\'e}aire. Il
reste {\`a} verifier que $\mathsf{Hec}_0$ poss{\`e}de des sommes directes. En
effet, montrons que la somme directe et le produit direct dans $\mathsf{Hec}_0$
sont les m{\^e}mes que dans $\dgcat$. Puisque l'on a 
$$ \mathsf{rep}_{ec}(\mathcal{A}\amalg \mathcal{B}, \mathcal{C}) \backsimeq
\mathsf{rep}_{ec}(\mathcal{A}, \mathcal{C})\times  \mathsf{rep}_{ec}(\mathcal{B}, \mathcal{C})$$
et que le groupe de Grothendieck pr{\'e}serve les produits, on a 
$$\mathrm{Hom}_{\mathsf{Hec}_0}(\mathcal{A}\amalg\mathcal{B},
  \mathcal{C}) \backsimeq \mathrm{Hom}_{\mathsf{Hec}_0}(\mathcal{A},
  \mathcal{C})\oplus \mathrm{Hom}_{\mathsf{Hec}_0}(\mathcal{B},
  \mathcal{C})\,.$$
Un argument analogue montre aussi que $\mathcal{A}\times\mathcal{B}$
  est bien le produit dans $\mathsf{Hec}_0$.
\end{proof}

\begin{lemme}\label{tens}
Le produit tensoriel $-\otimes^{\mathbb{L}}-$ de $\mathsf{Hec}$, voir remarque~\ref{monoidale}, induit
une structure mono{\"\i}dale sym{\'e}trique sur $\mathsf{Hec}_0$.
\end{lemme}
\begin{proof}
Soient $\mathcal{A}$, $\mathcal{B}$, $\mathcal{C}$ et $\mathcal{D}$
des dg-cat{\'e}gories. Comme dans la remarque~\ref{bitri}, le bifoncteur
$$
\xymatrix@!0 @R=1,5pc @C=12pc{
-\otimes^{\mathbb{L}}_k-: \mathsf{rep}_{ec}(\mathcal{A},
 \mathcal{B}) \times \mathsf{rep}_{ec}(\mathcal{C}, \mathcal{D})
\ar[r] & \mathsf{rep}_{ec}(\mathcal{A}\otimes\mathcal{C},
 \mathcal{B}\otimes \mathcal{D})\\
*+<2pc>{(X , Y)} \ar@{|->}[r]  & X_{cof}\otimes_k Y\, ,
}
$$
o{\`u} $X_{cof}$ est un remplacement cofibrant de $X$ dans la cat{\'e}gorie
$\mathrm{Mod}(\mathcal{A}^{op}\otimes^{\mathbb{L}}\mathcal{B})$,
par rapport {\`a} sa structure de cat{\'e}gorie de mod{\`e}les de Quillen, est
bi-triangul{\'e} puisqu'il provient aussi d'un dg-bifoncteur.
\end{proof}

\begin{remarque}
Soit $\mathsf{Heq}_{{ex}_0}$ la sous-cat{\'e}gorie pleine de
$\mathsf{Hec}_0$ dont les objets sont les dg-cat{\'e}gories
exactes. L'{\'e}quivalence entre $\mathsf{Hec}$ et $\mathsf{Heq}_{ex}$ permet d'{\'e}tablir une {\'e}quivalence entre $\mathsf{Hec}_0$ et $\mathsf{Heq}_{{ex}_0}$.
\end{remarque}

Soit $\mathcal{I}$ la dg-cat{\'e}gorie avec objets $1$, $2$ et dont les morphismes sont engendr{\'e}s par $m_{12}\in \mathrm{Hom}_{\mathcal{I}}^0(1, 2)$ soumis {\`a}
la relation $dm_{12}=0$. Pour une dg-cat{\'e}gorie
    $\mathcal{A}$, on note $\mbox{T}(\mathcal{A})$ la dg-cat{\'e}gorie
    $\mathcal{I}\otimes\mathcal{A}$. On dispose de deux inclusions
    canoniques $\mathcal{A} \stackrel{i_1}{\hookrightarrow} \mbox{T}(\mathcal{A})$ et $\mathcal{A} \stackrel{i_2}{\hookrightarrow} \mbox{T}(\mathcal{A})$ et d'une projection canonique
    $\mbox{T}(\mathcal{A})\rightarrow \mathcal{A}$.
\begin{remarque}\label{bimod}
La donn{\'e}e d'un $\mathcal{A}$-$\mathcal{B}$-bimodule $X$ est {\'e}quivalente
{\`a} la donn{\'e}e d'un dg-foncteur $X:\mathcal{A}\rightarrow
\mathrm{Mod}$-$\mathcal{B}$. Puisque la cat{\'e}gorie mono{\"\i}dale sym{\'e}trique
$(\dgcat,\otimes)$ est ferm{\'e}e, un
$\mathcal{A}$-$\mbox{T}(\mathcal{B})$-bimodule $X$ s'identifie {\`a} la
donn{\'e}e d'un morphisme de $\mathcal{A}$-$\mathcal{B}$-bimodules.
\end{remarque}
On note {\'e}galement $i_1$ et $i_2$ les morphismes de $\mathsf{Hec}$ associ{\'e}s
respectivement aux inclusions $i_1$ et $i_2$. Soit
$F:\mathsf{Hec}\rightarrow \mathsf{C}$ un foncteur {\`a} valeurs dans une
cat{\'e}gorie additive $\mathsf{C}$.

\begin{theoreme}\label{semiloc}
Les conditions suivantes sont {\'e}quivalentes~:
\begin{itemize}
\item[1)] Le foncteur $F$ est compos{\'e} d'un foncteur additif $\mathsf{Hec}_0\rightarrow \mathsf{C}$ et du foncteur canonique
  $\mathsf{Hec}\rightarrow \mathsf{Hec}_0$.
\item[2)] Pour toutes dg-cat{\'e}gories $\mathcal{A}$, $\mathcal{B}$, l'identit{\'e}
  $F(\left[X\right])+F(\left[Z\right])=F(\left[Y\right])$, est v{\'e}rifi{\'e}e
  dans $\mathrm{Hom}_{\mathsf{C}}(F(\mathcal{A}), F(\mathcal{B}))$ pour
  tout triangle $X \rightarrow Y \rightarrow Z \rightarrow X\left[1\right]$
  de $\mathsf{rep}_{ec}(\mathcal{A}, \mathcal{B})$.
\item[3)] Pour toute dg-cat{\'e}gorie $\mathcal{A}$, le morphisme 
$$
\xymatrix{
F(\mathcal{A})\oplus F(\mathcal{A})
\ar[rr]^-{\left[F(i_1)
        \, , \,F(i_2)\right]} && F(\mbox{T}(\mathcal{A}))
}
$$
est un isomorphisme dans $\mathsf{C}$.

\item[4)] Pout toute dg-cat{\'e}gorie pr{\'e}triangul{\'e}e $\mathcal{A}$ munie de sous-dg-cat{\'e}gories pleines pr{\'e}triangul{\'e}es
  $\mathcal{B}$ et $\mathcal{C}$ qui donnent lieu {\`a} une d{\'e}composition
  semi-orthogonale
  $\mathrm{H}^0(\mathcal{A})=(\mathrm{H}^0(\mathcal{B}),
  \mathrm{H}^0(\mathcal{C}))$, voir \cite{Bondal}, le morphisme
$$
F(\mathcal{B})\oplus F(\mathcal{C}) \rightarrow F(\mathcal{A})
$$
induit par les inclusions est un isomorphisme dans $\mathsf{C}$.

\end{itemize}
\end{theoreme}

\begin{proof}
Les conditions $1)$ et $2)$ sont {\'e}quivalentes par construction du groupe
de Grothendieck d'une cat{\'e}gorie triangul{\'e}e. Clairement la condition
$4)$ implique la condition $3)$. Montrons maintenant que la
condition $3)$ entra{\^\i}ne la condition $2)$. Pour une dg-cat{\'e}gorie
$\mathcal{E}$, notons $\mathrm{Mod}_{cf}\, \mathcal{E}$ la
sous-cat{\'e}gorie des objets cofibrants de
$\mathrm{Mod}\,\mathcal{E}$. Dans la suite, nous supposons que
$\mathcal{B}$ est cofibrant en tant que dg-cat{\'e}gorie. La cat{\'e}gorie
$\mathrm{Mod}\, \mbox{T}(\mathcal{B})$ s'identifie {\`a} la cat{\'e}gorie des
morphismes (ferm{\'e}s de degr{\'e} $0$)
$$ M_2 \stackrel{f}{\rightarrow} M_1$$
de $\mathrm{Mod}\, \mathcal{B}$. Les objets cofibrants correspondent
aux monomorphismes entre objets cofibrants de $\mathrm{Mod}\,
\mathcal{B}$. Pour un objet cofibrant $M$, on dispose d'une suite
exacte
$$ 0 \rightarrow M_2 \stackrel{f}{\rightarrow} M_1 \rightarrow
\mbox{Cok}f\rightarrow 0$$
de $\mathrm{Mod}\, \mathcal{B}$ fonctorielle en $M$. Notons
$$ 0 \rightarrow P_2 \rightarrow P_1 \rightarrow P_1/P_2 \rightarrow
0$$
la suite de dg-foncteurs de $\mathrm{Mod}_{cf}\,\mbox{T}(\mathcal{B})$ dans
$\mathrm{Mod}_{cf}\, \mathcal{B}$ obtenue ainsi. Clairement, les
foncteurs $P_2$, $P_1$, $P_1/P_2$ envoient les dg-modules
repr{\'e}sentables sur des dg-modules repr{\'e}sentables et donnent donc lieu
{\`a} des objets dans $\mathsf{rep}(\mbox{T}(\mathcal{B}), \mathcal{B})$. Soit
$$ X \rightarrow Y \rightarrow Z \rightarrow X\left[1\right]$$
un triangle de $\mathsf{rep}_{ec}(\mathcal{A}, \mathcal{B})$. Il est
isomorphe dans
$\mathcal{D}(\mathcal{A}^{op}\otimes^{\mathbb{L}}\mathcal{B})$ au
triangle associ{\'e} {\`a} une suite exacte courte
$$ 0 \rightarrow X' \stackrel{i}{\rightarrow} Y' \rightarrow Z'
\rightarrow 0$$
de bimodules cofibrants. Nous pouvons consid{\'e}rer le morphisme $M=(X'
\stackrel{i}{\rightarrow}Y')$ comme un dg-foncteur de
$\mathrm{Mod}_{cf}(\mathcal{A})$ dans $\mathrm{Mod}_{cf}\,
\mbox{T}(\mathcal{B})$. Clairement, il donne lieu {\`a} un objet de
$\mathsf{rep}(\mathcal{A}, \mbox{T}(\mathcal{B}))$. Nous avons
$$ P_2(M)=X', \, \, \, \, P_1(M)=Y', \, \,\, \, (P_1/P_2)(M)=Z'$$
dans $\mathsf{Hec}$. Pour montrer que l'on a $F(X')+F(Z')=F(Y')$, il
suffit donc de v{\'e}rifier que 
$$ F(P_2)+F(P_1/P_2)=F(P_1)\,.$$
Notons
$$
\xymatrix@!0 @R=1,5pc @C=5pc{
I_1: L \ar@{|->}[r] & (0\rightarrow L)\\
I_2:L \ar@{|->}[r] & ( L \stackrel{\mathbf{1}}{\rightarrow}L)
}
$$
 les dg-foncteurs de $\mathrm{Mod}_{cf}\, \mathcal{A}$ dans
 $\mathrm{Mod}_{cf}\, \mbox{T}(\mathcal{A})$ induits par $i_1$ et
 $i_2$. Clairement, nous avons
$$
(P_1/P_2) \circ I_1=\mathbf{1}, \, \, \, (P_1/P_2) \circ I_2=0, \, \, \,
P_2 \circ I_1=0, \, \, \, P_2 \circ I_2=\mathbf{1}$$
dans $\mathsf{Hec}$.
Il s'ensuit que dans la cat{\'e}gorie additive $\mathsf{C}$, nous avons
$$
\begin{bmatrix}F(P_1/P_2)\\
F(P_2)
\end{bmatrix}
\circ
\begin{bmatrix}
F(I_1) & F(I_2)
\end{bmatrix}
=
\begin{bmatrix}
\mathbf{1} & 0 \\
0 & \mathbf{1}
\end{bmatrix}
$$
et donc
$$
(F(P_1/P_2)+F(P_2))\circ \begin{bmatrix}F(I_1) & F(I_2) \end{bmatrix}
= \begin{bmatrix}\mathbf{1} & \mathbf{1}\end{bmatrix}\,.$$
De l'autre c{\^o}t{\'e}, nous avons
$$P_1 \circ I_1=\mathbf{1}, \, \, \, P_1 \circ I_2=\mathbf{1}$$
dans $\mathsf{Hec}$ et donc
$$
F(P_1)\circ \begin{bmatrix}F(I_1) & F(I_2)
\end{bmatrix}=\begin{bmatrix}\mathbf{1} &  \mathbf{1}\end{bmatrix}$$ 
dans $\mathsf{C}$. Comme $ \begin{bmatrix} F(I_1) & F(I_2)\end{bmatrix}$ est
inversible, il s'ensuit que l'on a bien
$$
F(P_1/P_2)+F(P_2)=F(P_1)\,.$$

Montrons maintenant que la condition $1)$ implique la condition
$3)$. Ecrivons $F$ comme compos{\'e} d'un foncteur additif
$\overline{F}:\mathsf{Hec}_0 \rightarrow \mathsf{C}$ et du foncteur
canonique $\mathsf{Hec}\rightarrow \mathsf{Hec}_0$. Nous avons
$$
\overline{F}\,\mathcal{A}\oplus \overline{F}\,\mathcal{A}\simeq
\overline{F}(\mathcal{A}\oplus \mathcal{A})\, .
$$
Il suffit donc de montrer que dans $\mathsf{Hec}_0$, le morphisme
canonique
$$
\mathcal{A}\oplus\mathcal{A} \longrightarrow \mbox{T}(\mathcal{A})
$$
devient inversible. Par le lemme de Yoneda, il suffit de montrer que
pour toute dg-cat{\'e}gorie $\mathcal{U}$, l'application
$$
\mbox{K}_0(\mathsf{rep}_{ec}(\mathcal{U},\mathcal{A}))\oplus
\mbox{K}_0(\mathsf{rep}_{ec}(\mathcal{U},\mathcal{A})) \longrightarrow
\mbox{K}_0(\mathsf{rep}_{ec}(\mathcal{U}, \mbox{T}(\mathcal{A})))
$$
est bijective. Ceci r{\'e}sulte du fait que la cat{\'e}gorie triangul{\'e}e
$\mathsf{rep}_{ec}(\mathcal{U}, \mbox{T}(\mathcal{A}))$ admet une
d{\'e}composition semi-orthogonale en deux sous-cat{\'e}gories {\'e}quivalentes {\`a}
$\mathsf{rep}_{ec}(\mathcal{U}, \mathcal{A})$. En effet, si
$X$ est un objet de $\mathsf{rep}_{ec}(\mathcal{U},
\mbox{T}(\mathcal{A}))$ qui est cofibrant dans 
$$
\mathrm{Mod}(\mathcal{U}^{op}\otimes
\mbox{T}(\mathcal{A}))\stackrel{\sim}{\longrightarrow}
\mathrm{Mod}(\mbox{T}(\mathcal{U}^{op}\otimes \mathcal{A}))\, ,$$
nous pouvons l'identifier avec un monomorphisme $X_2
\stackrel{i}{\rightarrow} X_1$ entre objets cofibrants appartenant {\`a}
$\mathsf{rep}_{ec}(\mathcal{U},\mathcal{A})$. Alors le triangle
associ{\'e} {\`a} $X$ par la d{\'e}composition semi-orthogonale est induit par la
suite exacte de morphismes
$$
\xymatrix{
0 \ar[r] \ar[d] & X_2 \ar@{=}[r] \ar@{=}[d] & X_2 \ar[d]^-{i} \ar[r] & 0
\ar[d] \ar[r] & 0 \ar[d]\\
0 \ar[r] & X_2 \ar[r]_-{i} & X_1 \ar[r] & X_1/X_2 \ar[r] & 0\,.
}
$$
Un argument analogue au pr{\'e}c{\'e}dent montre que la condition $1)$
implique la condition $4)$.
\end{proof}

\section{dg-foncteurs de Morita}

Soit $\mathcal{T}$ une cat{\'e}gorie triangul{\'e}e. On note
$\mathcal{T}^{\,\kar}$ sa compl{\'e}tion idempotente, voir \cite{Balmer}.
Ici, l'exposant repr{\'e}sente la moiti{\'e} du symbole $\oplus$.
Un dg-foncteur $F$ de $\mathcal{C}$ vers $\mathcal{E}$ est {\em de Morita} si:
\begin{itemize}
\item[-] pour tous objets $c_1$ et $c_2$ dans $\mathcal{C}$, le morphisme
  de complexes de $\mathrm{Hom}_{\mathcal{C}}(c_1, c_2)$ vers
  $\mathrm{Hom}_{\mathcal{E}}(F(c_1),F(c_2))$ est un
  quasi-isomorphisme et
\item[-] le foncteur $\mathrm{H}^0(\mbox{pre-tr}(F))^{\kar}$ de
  $\mathrm{H}^0(\mbox{pre-tr}(\mathcal{C}))^{\kar}$ vers
  $\mathrm{H}^0(\mbox{pre-tr}(\mathcal{E}))^{\kar}$ est essentiellement surjectif. 
\end{itemize}

\begin{remarque}
Cette notion est importante car un dg-foncteur $F:\mathcal{C}
\rightarrow \mathcal{E}$ est de Morita si et seulement s'il induit une
{\'e}quivalence
$$F_*:\mathcal{D}(\mathcal{E})
\stackrel{\sim}{\rightarrow}\mathcal{D}(\mathcal{C})$$
dans les cat{\'e}gories d{\'e}riv{\'e}es, voir \cite{Keller}.
\end{remarque}

On introduira une structure de cat{\'e}gorie de mod{\`e}les de
Quillen {\`a} engendrement cofibrant dans $\dgcat$ dont les
{\'e}quivalences faibles sont les dg-foncteurs de Morita. Pour cela, on se servira du th{\'e}or{\`e}me~2.1.19
de \cite{Hovey}.
Soit $\mathrm{idemh}_0$ la dg-cat{\'e}gorie avec un seul object $0$ et
dont les morphismes sont engendr{\'e}s par des endomorphismes $e$ de degr{\'e}
$0$ et $h$ de degr{\'e} $-1$ soumis aux relations $de=0$
et $dh=e^2-e$.
Pour un entier $n>0$ et des entiers $k_0, \ldots , k_n$, soit
$\mathrm{idemh}_n(k_0, \ldots , k_n)$ la dg-cat{\'e}gorie obtenue {\`a} partir
de la dg-cat{\'e}gorie $\mathcal{M}_n(k_0, \ldots , k_n)$ en rajoutant de
nouveaux generateurs et relations~: on rajoute, pour tous $0 \leq i, j
\leq n$, des morphismes $e_{i, j}$ de source $j$ et de but $i$ de
degr{\'e} $k_i-k_j$ soumis aux relations
$$d(E)=0 $$
o{\`u} $E$ est la matrice de coefficients $e_{i, j}$.
On rajoute aussi, pour tous $0 \leq i, j \leq n$, des morphismes $h_{i,
  j}$ de source $j$ et de but $i$ de degr{\'e} $k_i-k_j-1$ soumis aux
relations
$$d(H)=E^2-E$$
o{\`u} $H$ est la matrice de coefficients $h_{i, j}$.
Pour $n\geq0$, soit $\mathrm{Facth}_n(k_0, \ldots , k_n)$ la sous-dg-cat{\'e}gorie pleine
de la dg-cat{\'e}gorie des dg-modules ({\`a} droite) sur
$\mathrm{idemh}_n(k_0, \ldots , k_n)$ dont les objets sont les
$\hat{l}$, $0 \leq l \leq n$, le c{\^o}ne it{\'e}r{\'e} $X_n$, qui a le m{\^e}me
module gradu{\'e} sous-jacent que le dg-module 
$$ X=\bigoplus^n_{l=0} \hat{l}\left[ k_l \right]$$
et dont la diff{\'e}rentielle est $d_X+\widehat{Q}$, le dg-module
$\bigoplus_{i \in \mathbb{N}}Y_i$, o{\`u} $Y_i=X_n$ pour tout $i \in \mathbb{N}$, et le c{\^o}ne $C_n$ sur le morphisme
$$ \Phi:\bigoplus_{i \in \mathbb{N}}Y_i \longrightarrow
\bigoplus_{i \in \mathbb{N}}Y_i$$ dont les composantes sont toutes de la forme
$$
\xymatrix{
Y_j \ar[rr]^-{\left[Id-E \,\, E \right]^t} && Y_j \oplus Y_{j+1}
    \ar@{^{(}->}[rr]^-{can} && \bigoplus_{i \in \mathbb{N}}Y_i\,.
}$$
\begin{remarque} \label{black-box}
La classe $\left[\Phi \right]$ du morphisme $\Phi$ dans
$\mathrm{H}^0(\mathrm{Facth}_n(k_0, \ldots , k_n))$ admet une
r{\'e}traction $\left[\Psi \right]$, o{\`u} la premi{\`e}re composante de $\Psi$
est $Id-E$ et les suivantes sont toutes de la forme
$$\xymatrix{
Y_j \ar[rr]^-{\left[ E \,\, Id-E \right]^t} &&  Y_j \oplus Y_{j+1}
    \ar@{^{(}->}[rr]^-{can} && \bigoplus_{i \in \mathbb{N}} Y_i\,.
}$$
Cela implique que dans la cat{\'e}gorie homotopique des dg-modules ({\`a}
droite) sur $\mathrm{idemh}_n(k_0, \ldots , k_n)$, le morphisme
$\left[E \right]$ admet une image et cette image est isomorphe {\`a} $C_n$. 
\end{remarque}

Soit $L:\mathcal{A} \rightarrow \mathrm{Facth}_n(k_0, \ldots , k_n)$ le dg
foncteur qui envoie $3$ sur $C_n$. On consid{\`e}re la somme amalgam{\'e}e
$$\xymatrix{
\mathcal{A} \ar[r]^-{L} \ar[d]_F \ar@{}[dr]|{\lrcorner} & \mathrm{Facth}_n(k_0, \ldots , k_n) \ar[d] \\
\mathcal{K} \ar[r] & \mathrm{Facth}_n(k_0, \ldots ,
k_n) \amalg_{\mathcal{A}} \mathcal{K} 
}$$
et on d{\'e}finit $\mathrm{facth}_n(k_0, \ldots , k_n)$ comme la sous-dg-cat{\'e}gorie
pleine de la somme amalgam{\'e}e dont les objets sont les images des
objets $\hat{l}$, $0 \leq l \leq n$, de $\mathrm{Facth}_n(k_0, \ldots , k_n)$
et l'image $C\!h_{\,n}$ de l'object $2$ de $\mathcal{K}$. On note
$L_n(k_0, \ldots , k_n)$ le dg foncteur fid{\`e}le, mais non plein, de 
$\mathrm{idemh}_n(k_0, \ldots , k_n)$ dans $\mathrm{facth}_n(k_0, \ldots , k_n)$.

\begin{theoreme}\label{theorem2}
 Si on consid{\`e}re pour cat{\'e}gorie $\mathcal{C}$ la cat{\'e}gorie $\dgcat$,
pour classe $W$ la sous-cat{\'e}gorie de $\dgcat$ des dg-foncteurs
de Morita, pour classe $J$ les dg-foncteurs $F$ et $R(n)$,
$n\in \mathbb{Z}$, (voir \cite{cras}), $F(n)$, $n \in \mathbb{Z}$, $I_n(k_0,
 \ldots , k_n)$ et $L_n(k_0, \ldots , k_n)$ et pour classe $I$ les dg-foncteurs $Q$ et $S(n)$,
$n\in \mathbb{Z}$, (voir \cite{cras}), alors les conditions du th{\'e}or{\`e}me \cite[2.1.19]{Hovey} sont satisfaites.
\end{theoreme}

\begin{remarque}
On peut montrer ais{\'e}ment que pour la structure obtenue les objets
fibrants sont les dg-cat{\'e}gories $\mathcal{A}$ telles que l'image
essentielle du plongement $\mathrm{H}^0(\mathcal{A}) \hookrightarrow
\mathcal{D}(\mathcal{A})$ est stable par suspension, c{\^o}nes et facteurs
directs.
\end{remarque}

On observe facilement que les conditions $\mbox{{\it (i)}}$,
$\mbox{{\it (ii)}}$ et $\mbox{{\it (iii)}}$ du th{\'e}or{\`e}me~2.1.19
de \cite{Hovey} sont verifi{\'e}es.

\begin{lemme}\label{J-inj-Surj2}
$\mbox{{\bf Surj}} = J-\mbox{inj} \cap W\,.$
\end{lemme}
\begin{proof}
Montrons l'inclusion $\supseteq$.
Soit $H$ un dg-foncteur de $\mathcal{D}$ vers
$\mathcal{S}$ qui appartient {\`a} $J-\mbox{inj} \cap W$. La classe
$R(n)-\mbox{drt}$ est form{\'e}e des dg-foncteurs surjectifs au niveau des
complexes de morphismes. Comme $H \in W$, il suffit de montrer que $H$
est surjectif au niveau des objets. Soit $A \in \mathcal{S}$ un
object quelconque. Comme $H \in W$, il existe $C \in
\mathrm{H}^0(\mbox{pre-tr}(\mathcal{D}))^{\kar}$ et un
isomorphisme $q$ de
$\mathrm{H}^0(\mbox{pre-tr}(\mathcal{S}))^{\kar}$ entre l'image
par $\mathrm{H}^0(\mbox{pre-tr}(\mathcal{H}))^{\kar}$ de $C$ et
$A$. On identifie la dg-cat{\'e}gorie
$\mathrm{H}^0(\mbox{pre-tr}(\mathcal{D}))$ avec son image dans
$\mathrm{H}^0(\mbox{pre-tr}(\mathcal{D}))^{\kar}$. Il existe
alors deux possibilit{\'e}s~:
\begin{enumerate}
\item l'object $C$ appartient {\`a} la dg-cat{\'e}gorie
  $\mathrm{H}^0(\mbox{pre-tr}(\mathcal{D}))$. Alors on est dans les
  conditions du th{\'e}or{\`e}me pr{\'e}c{\'e}dent.
\item L'object $C$ appartient {\`a} la dg-cat{\'e}gorie
  $\mathrm{H}^0(\mbox{pre-tr}(\mathcal{D}))^{\kar}$ mais non pas {\`a} la
  cat{\'e}gorie $\mathrm{H}^0(\mbox{pre-tr}(\mathcal{D}))$. Alors $C$ est de la forme
  $C=(\overline{C}, \left[ \psi \right])$, o{\`u} $ \overline{C} \in
  \mbox{pre-tr}(\mathcal{D})$ et $\psi$ est un endomorphisme ferm{\'e} de
  $C$ dans $\mbox{pre-tr}(\mathcal{D})$ tel que $\left[
  \psi^2\right] = \left[\psi \right]$. Cette information permet de construire le
  diagramme commutatif suivant
$$\xymatrix @R=1pc{
*+<1pc>{\mathrm{idemh}_n(k_0, \ldots , k_n)} \ar@{^{(}->}[dd]_{can} \ar[r]
& \mathcal{D} \ar[d]^H \\
 &  \mathcal{S} \ar[d] \\
\mathrm{Facth}_n(k_0, \ldots , k_n) \ar[r]^-T & \mathrm{Mod}\,\mathcal{S}\,.
}$$
Gr{\^a}ce {\`a} la remarque~\ref{black-box}, l'isomorphisme $q$ dans $\mathrm{H}^0(\mbox{pre-tr}(\mathcal{D}))^{\kar}$ permet de
construire un dg-foncteur de $\mathcal{K}$ vers
$\mathrm{Mod}\,\mathcal{S}$ qui envoie $1$ sur
$T(C\!h_{\,n})$ et $2$ sur $A$. On peut donc {\'e}tendre le
dg-foncteur $T$ {\`a} un dg-foncteur $\overline{T}$ de $\mathrm{Facth}_n(k_0, \ldots ,
k_n) \amalg_{\mathcal{A}} \mathcal{K}$ vers
$\mathrm{Mod}\,\mathcal{S}$, tel que l'image de sa restriction {\`a} la
sous-dg-cat{\'e}gorie pleine $\mathrm{facth}_n(k_0, \ldots , k_n)$, qu'on note $\widetilde{T}$, est dans $\mathcal{S}$. On dispose donc
du diagramme commutatif suivant  
$$\xymatrix@!0 @R=3pc @C=4pc{
 & & \mathrm{idemh}_n \ar[rr] \ar[dr]_-{L_n} \ar[dl]^-{can} & 
 & \mathcal{D} \ar[d]^(0.3){H} \\
\mathcal{A} \ar@{}[drr]|{\lrcorner} \ar[r] \ar[dr]_-F & \mathrm{Facth}_n \ar[dr] \ar[drrr]^(0.3){T} & &
 \mathrm{facth}_n \ar[r]^-{\widetilde{T}}
 \ar@{^{(}->}[dl]|\hole \ar@{.>}[ur]  & \mathcal{S} \ar[d] \\
 & \mathcal{K} \ar[r] & \mathrm{Facth}_n \amalg_{\mathcal{A}} \mathcal{K} \ar[rr]_-{\overline{T}} &  & \mbox{Mod}\,\mathcal{S}\,.
}$$

\end{enumerate}

Montrons maintenant l'inclusion $\subseteq$. Soit $H$ un dg-foncteur
de $\mathcal{N}$ vers $\mathcal{E}$ dans la classe $\mbox{{\bf Surj}}$. Comme $H$ est surjectif au niveau des objets et un
quasi-isomorphisme au niveau des complexes de morphismes, on a $H \in
W$. La classe $R(n)-\mbox{drt}$ est form{\'e}e des dg-foncteurs surjectifs
au niveau des complexes de morphismes. Il suffit de montrer que $H$ a
la propri{\'e}t{\'e} de rel{\`e}vement {\`a} droite par rapport {\`a} $F$, $F(n)$, $n \in
\mathbb{Z}$, $I_n(k_0, \ldots , k_n)$ et $L_n(k_0, \ldots , k_n)$. On
consid{\`e}re ces deux cas~:
\begin{enumerate}
\item on sait par le th{\'e}or{\`e}me pr{\'e}c{\'e}dent que  $H \in F-\mbox{drt}$,
  $H \in F(n)-\mbox{drt}$ et $H \in I_n(k_0, \ldots , k_n)-\mbox{drt}$.  
\item On consid{\`e}re le diagramme suivant~:
$$\xymatrix{
\mathrm{idemh}_n(k_0,\ldots, k_n) \ar[r] \ar[d]_{L_n(k_0,
  \ldots , k_n)} & \mathcal{N} \ar[d]^H\\
\mathrm{facth}_n(k_0, \ldots , k_n) \ar[r]_-{T} & \mathcal{E}
}$$
Le lemme~\ref{lemme-clef} permet d'{\'e}tendre le dg-foncteur $T$ en un dg-foncteur
$\overline{T}$ de $\mathrm{Facth}_n(k_0, \ldots ,
 k_n)\amalg_{\mathcal{A}} \mathcal{K}$ vers
 $\mathrm{Mod}\,\mathcal{E}$. Puisque $H$ appartient {\`a} la classe
 $\mbox{{\bf Surj}}$, on peut relever le dg-foncteur $\overline{T}$
 vers $\mbox{Mod}\,\mathcal{N}$. On remarque que la restriction de ce
 rel{\`e}vement {\`a} la sous-dg-cat{\'e}gorie pleine $\mathrm{facth}_n(k_0,
 \ldots , k_n)$ fournit un rel{\`e}vement du dg-foncteur $T$ vers la
 dg-cat{\'e}gorie $\mathcal{N}$. On dispose du diagramme commutatif
 suivant
$$\xymatrix@!0 @R=2pc @C=4pc{
&  & \mathrm{idemh}_n \ar[rr] \ar[dl] \ar[dd]|\hole & &  *+<1pc>{\mathcal{N}}
\ar[dd]^-H \ar@{^{(}->}[dl] \\
\mathcal{A} \ar@{}[ddr]|{\lrcorner} \ar[dd]_F \ar[r] & \mathrm{Facth}_n \ar[dd] \ar[rr]  & &
\mbox{Mod}\,\mathcal{N} \ar[dd] & \\
& & \mathrm{facth}_n \ar[rr]|\hole^(0.3){T} \ar[dl] &    & *+<1pc>{\mathcal{E}}
\ar@{^{(}->}[dl] \\
\mathcal{K} \ar[r] & \mathrm{Facth}_n \amalg_{\mathcal{A}} \mathcal{K}
\ar[rr]^-{\overline{T}} &   & \mbox{Mod}\,\mathcal{E} & 
}$$
\end{enumerate}
\end{proof}

\begin{lemme}\label{J-cell-dans-W2}
 On a $J-\mbox{cell}\subset W$.
\end{lemme}
\begin{proof}
On sait d{\'e}j{\`a} par le th{\'e}or{\`e}me pr{\'e}c{\'e}dent que les classes
$F-\mbox{cell}$, $R(n)-\mbox{cell}$, $F(n)-\mbox{cell}$ et $I_n(k_0,
\ldots , k_n)-\mbox{cell}$ sont form{\'e}es de dg-foncteurs quasi-{\'e}quiconiques donc contenues dans la classe $W$.
 
Soit maintenant $T:\mathrm{idemh}_n(k_0,\ldots, k_n) \rightarrow \mathcal{J}$ un
dg-foncteur quelconque. On consid{\`e}re la somme amalgam{\'e}e suivante
$$\xymatrix{
\mathrm{idemh}_n(k_0,\ldots, k_n) \ar[d]_{L_n(k_0, \ldots ,
  k_n)} \ar[r]^-{T}  \ar@{}[dr]|{\lrcorner} & \mathcal{J} \ar[d]^{inc} \\
\mathrm{facth}_n(k_0, \ldots , k_n) \ar[r] & \mathcal{U}
}$$ 
dans $\dgcat$. Il s'agit de v{\'e}rifier que $inc \in W$. Il
faut verifier deux conditions~:
\begin{enumerate}
\item L'inclusion $\mathrm{Hom}_{\mathcal{J}}(X,Y) \rightarrow \mathrm{Hom}_{\mathcal{U}} (inc(X), inc(Y))$ est un quasi-isomorphisme
  pour tous $X,Y \in \mathcal{J}$. On dispose du diagramme commutatif
  suivant associ{\'e} {\`a} la construction de la somme amalgam{\'e}e.
$$\xymatrix@!0 @R=2pc @C=4pc{
&  & \mathrm{idemh}_n \ar[rr]^T \ar[dl] \ar[dd]|\hole & &  \mathcal{J}
\ar[dd]^-{inc} \ar[dl]^L \\
\mathcal{A} \ar@{}[ddr]|{\lrcorner} \ar[dd]_F \ar[r] & \mathrm{Facth}_n \ar[dd] \ar[rr]  & & \mathcal{N} \ar[dd]^(0.3){H} & \\
& & \mathrm{facth}_n \ar[rr]|\hole \ar@{^{(}->}[dl] &    & *+<1pc>{\mathcal{U}}
\ar@{^{(}->}[dl] \\
\mathcal{K} \ar[r] & \mathrm{Facth}_n \amalg_{\mathcal{A}} \mathcal{K}
\ar[rr] &   & \mathcal{E} & 
}$$
 En effet, on remarque que le dg-foncteur $inc$ s'identifie {\`a} la
 composition de $L$ et $H$ suivi d'une restriction. On peut identifier
 la dg-cat{\'e}gorie $\mathcal{N}$ {\`a} la sous-dg-cat{\'e}gorie pleine des
 dg-modules ({\`a} droite) sur $\mathcal{J}$ dont les objets sont les
 dg-modules represent{\'e}es, le $n$-i{\`e}me c{\^o}ne it{\'e}r{\'e} associ{\'e} {\`a} l'image par
 le dg-foncteur $T$, la somme $\bigoplus_{\mathbb{N}}T(X_n)$ et
 $T(C\!h_{\,n})$. Le dg-foncteur $L$ s'identifie donc au plongement
 de Yoneda et est donc pleinement fid{\`e}le. Le dg-foncteur $H$ s'identifie au dg-foncteur de $\mathcal{N}$ vers
 la somme amalgam{\'e} avec $\mathcal{K}$. On sait d'apr{\`e}s \cite{cras}
 que $H$ est donc un quasi-isomorphisme. Cela implique que $inc$
 satisfait la condition.
\item  Le foncteur
  $\mathrm{H}^0(\mbox{pre-tr}(inc))^{\kar}:\mathrm{H}^0(\mbox{pre-tr}(\mathcal{J}))^{\kar}
  \rightarrow \mathrm{H}^0(\mbox{pre-tr}(\mathcal{U})))^{\kar}$ est essentiellement
  surjectif. En effet, la dg-cat{\'e}gorie $\mathcal{U}$ poss{\`e}de un object
  $C\!h_{\,n}$ de plus que la dg-cat{\'e}gorie $\mathcal{J}$. Or l'objet
  $C\!h_{\,n}$ est isomorphe dans
  $\mathrm{H}^0(\mbox{pre-tr}(\mathcal{U})^{\kar}$ {\`a} l'image par le
  dg-foncteur $inc$ de l'object $T(X_n)$ muni de l'endomorphisme $T(E)$.  
\end{enumerate}
\end{proof}
Nous avons v{\'e}rifi{\'e} que $J-\mbox{cell}\subseteq W$ (lemme~\ref{J-cell-dans-W2}) et
que $I-\mbox{inj}$ est {\'e}gal {\`a} $J-\mbox{inj}\cap W$ (lemmes~\ref{I-inj-Surj} et \ref{J-inj-Surj2}).
Ces conditions impliquent celles du th{\'e}or{\`e}me de Hovey \cite[2.1.19]{Hovey}.

On note $\mathcal{M}or$ la classe des dg-foncteurs de Morita et
$\mathsf{Hmo}$ la cat{\'e}gorie homotopique de $\dgcat$ par raport {\`a}
cette classe. Soient $\mathcal{A}$ et $\mathcal{B}$ des dg-cat{\'e}gories et $\mbox{can}_2: \mathsf{Heq} \rightarrow \mathsf{Hmo}$ le
foncteur canonique.

\begin{corollary} \label{adjonct2}
On a une adjonction
$$\mathrm{Hom}_{\mathsf{Hmo}}(\mbox{can}_2(\mathcal{A}),
\mathcal{B}) \stackrel{\sim}{\longrightarrow} \mathrm{Hom}_{\mathsf{Heq}}(\mathcal{A},
\mathcal{B}_{fib})\,,
$$
o{\`u} $\mathcal{B}_{fib}$ est un remplacement fibrant
de la dg-cat{\'e}gorie $\mathcal{B}$ par rapport {\`a} la structure de cat{\'e}gorie de mod{\`e}les de Quillen d{\'e}crite
dans le th{\'e}or{\`e}me~\ref{theorem2}.
\end{corollary}
\begin{proof}
En effet, les structures de cat{\'e}gorie de mod{\`e}les de Quillen du th{\'e}or{\`e}me~\ref{theorem} et du th{\'e}or{\`e}me~\ref{theorem2} ont les m{\^e}mes
cofibrations et les m{\^e}mes fibrations acycliques. Cela implique que les
objets cofibrants et les objets cylindres sont les m{\^e}mes et par
cons{\'e}quent la relation d'homotopie est la m{\^e}me dans les deux
situations.
On a donc des bijections
$$
\mathrm{Hom}_{\mathsf{Hmo}}(\mbox{can}_2(\mathcal{A}), \mathcal{B})
\stackrel{\sim}{\longrightarrow} \mathrm{Hom}_{\dgcat}(\mathcal{A}_{cof},\mathcal{B}_{fib})/ \mbox{htp} \stackrel{\sim}{\longleftarrow} \mathrm{Hom}_{\mathsf{Heq}}(\mathcal{A}, \mathcal{B}_{fib})\,.
$$
\end{proof}
Soit $\mathsf{Heq}_{fib}$ la sous-cat{\'e}gorie pleine de $\mathsf{Heq}$
form{\'e}e des objets fibrants par rapport {\`a} la structure de cat{\'e}gorie de
mod{\`e}les de Quillen d{\'e}crite dans le th{\'e}or{\`e}me~\ref{theorem2}. On note
$\mbox{inc}:\mathsf{Heq}_{fib}\hookrightarrow \mathsf{Heq}$ l'inclusion
{\'e}vidente. Le corollaire~\ref{adjonct2} a comme consequence le

\begin{corollary}\label{adjonct5}
On a une adjonction
$$
\mathrm{Hom}_{\mathsf{Heq}}(\mathcal{A}, \mbox{inc}(\mathcal{B})) \stackrel{\sim}{\longleftarrow}
\mathrm{Hom}_{\mathsf{Heq}_{fib}}(\mathcal{A}_{fib}, \mathcal{B})\,.
$$
\end{corollary}

Les corollaires~\ref{adjonct2} et~\ref{adjonct5} {\'e}tablissent une
equivalence de cat{\'e}gories entre $\mathsf{Hmo}$ et
$\mathsf{Heq}_{fib}$.
On a donc le diagramme
$$
\xymatrix@!0 @R=4pc @C=8pc{
*+<0.5pc>{\mathsf{Heq}} \ar@<-1ex>[d]_-{can_2} \ar@<-1ex>[r]_-{(-)_{fib}} &
*+<1pc>{\mathsf{Heq}_{fib}} 
\ar@<2.5ex>[dl]^-{\sim} \ar@<-1ex>@{^{(}->}[l]_-{inc}\\
*+<0.5pc>{\mathsf{Hmo}} \ar@<-1ex>[u]_-{(-)_{fib}} \ar@<-1.5ex>[ur]^-{\sim} & 
}
$$

Soit $\mathsf{Heq}_{kar}$ la sous-cat{\'e}gorie pleine de
$\mathsf{Heq}_{ex}$  form{\'e}e des dg-cat{\'e}gories $\mathcal{A}$ telles
que l'image du plongement $h:\mathrm{H}^0(\mathcal{A})
\hookrightarrow \mathcal{D}(\mathcal{A})$ est stable par facteurs
directs.

\begin{lemme}\label{lemme3}
Les cat{\'e}gories $\mathsf{Heq}_{kar}$ et $\mathsf{Heq}_{fib}$ sont {\'e}quivalentes.
\end{lemme}
\begin{proof}
C'est une cons{\'e}quence de la remarque~\ref{fibrants}, qui implique
que le plongement de $\mathsf{Heq}_{ex}$ dans la sous-cat{\'e}gorie pleine
de $\mathsf{Heq}$ form{\'e}e des objets fibrants par rapport {\`a} la structure de cat{\'e}gorie de
mod{\`e}les de Quillen d{\'e}crite dans le th{\'e}or{\`e}me~\ref{theorem} est une
{\'e}quivalence de cat{\'e}gories.
\end{proof}

On note $\mathsf{rep}_{mor}(\mathcal{A}, \mathcal{B})$ la sous-cat{\'e}gorie pleine
de la cat{\'e}gorie d{\'e}riv{\'e} $\mathcal{D}(\mathcal{A}^{op} \otimes^{\mathbb{L}} \mathcal{B})$,
  voir \cite{Keller} \cite{Drinfeld}, dont les objets sont les dg
  $\mathcal{A}\mbox{-}\mathcal{B}\,$-bimodules $X$ telles que $X(?, A)$ est
  isomorphe dans $\mathcal{D}(\mathcal{B})$ {\`a} un objet de l'image du
  dg-foncteur canonique $\mathrm{H}^0(\mbox{pre-tr}(\mathcal{B}))^{\kar}
  \rightarrow
  \mathrm{H}^0(\mathrm{Mod}\, \mathcal{B})$, pour tout $A \in \mathcal{A}$.  

\begin{corollary}\label{adjonct6}
On a une bijection
$$\mathrm{Hom}_{\mathsf{Hmo}}(\mathcal{A},
\mathcal{B}) \stackrel{\sim}{\longrightarrow} \mathrm{Iso}(\mathsf{rep}_{mor}(\mathcal{A}, \mathcal{B}))\,.
$$
\end{corollary}

\begin{proof}
D'apr{\`e}s le corollaire~\ref{adjonct2}, $\mathrm{Hom}_{\mathsf{Hmo}}(\mathcal{A},
\mathcal{B})$ s'identifie {\`a}  $\mathrm{Hom}_{\mathsf{Heq}}(\mathcal{A},
\mathcal{B}_{fib})$. Par \cite{Toen} on sait que
ce dernier ensemble s'identifie {\`a}
$\mathrm{Iso}(\mathsf{rep}(\mathcal{A},\mathcal{B}_{fib}))$.
On sait par le lemme~\ref{lemme3} que $\mathcal{B}_{fib}$ est
isomorphe dans $\mathsf{Heq}$ {\`a} une dg-cat{\'e}gorie qui appartient {\`a}
$\mathsf{Heq}_{kar}$. Cela implique que $\mathsf{rep}(\mathcal{A},\mathcal{B}_{fib})$
s'identifie {\`a} $\mathsf{rep}_{mor}(\mathcal{A},\mathcal{B})$, d'o{\`u} le resultat.
\end{proof}

\begin{remarque}\label{monoi}
On sait que les structures de cat{\'e}gorie de mod{\`e}les de Quillen dans
$\dgcat$ de \cite{cras} et du th{\'e}or{\`e}me~\ref{theorem2} ont les m{\^e}mes
cofibrations et les m{\^e}mes fibrations acycliques. Cela a pour
consequence que
\begin{itemize}
\item[-] Le foncteur de remplacement cofibrant simplicial $\Gamma^{\ast}$, voir
\cite{Hovey}, dans $\dgcat$ est le m{\^e}me dans les deux situations. Cela
implique que les corollaires~\ref{adjonct2}, ~\ref{adjonct5} et~\ref{adjonct6} sont encore vrais si on remplace $\mathrm{Hom}$ par
l'espace de morphismes, $\underline{\mathrm{Map}}$, voir \cite{Hovey}, et
$\mathrm{Iso}(\mathsf{rep}_{mor}(\mathcal{A}, \mathcal{B}))$ dans le
corollaire~\ref{adjonct6}, par le nerf de la sous-cat{\'e}gorie de
$\mathrm{Mod}(\mathcal{A}^{op}\otimes^{\mathbb{L}}\mathcal{B})$ dont les objets
sont les m{\^e}mes que ceux de $\mathsf{rep}_{mor}(\mathcal{A},
\mathcal{B})$ et dont les morphismes sont les quasi-isomorphismes.
\item[-] La cat{\'e}gorie mono{\"\i}dale sym{\'e}trique
  $(\mathsf{Hmo},\,-\otimes^{\mathbb{L}}-)$ est bien d{\'e}finie et ferm{\'e}e, voir \cite{Hovey}, et l'espace de
  morphismes interne de $\mathsf{Hmo}$, $\mathbb{R}\underline{\mathrm{Hom}}_{\mathsf{Hmo}}(\mathcal{A},
    \mathcal{B})$ s'identifie {\`a}
  $\mathbb{R}\underline{\mathrm{Hom}}_{\mathsf{Heq}}(\mathcal{A}, \mathcal{B}_{fib})$, voir \cite{Toen}.
\end{itemize} 
\end{remarque}

\section{Invariants}
Soient $\mathcal{A}$, $\mathcal{B}$ et $\mathcal{C}$ des
dg-cat{\'e}gories. Soit $\mathsf{Hmo}_0$ la cat{\'e}gorie qui a pour objets
les petites dg-cat{\'e}gories et telle que $\mbox{Hom}_{\mathsf{Hmo}_0}(\mathcal{A}, \mathcal{B})$ est le groupe
de Grothendieck de la cat{\'e}gorie triangul{\'e}e
$\mathsf{rep}_{mor}(\mathcal{A}, \mathcal{B})$. Un argument
compl{\`e}tement analogue {\`a} la remarque~\ref{bitri} montre que l'op{\'e}ration
de composition dans $\mathsf{Hmo}_0$ est bien d{\'e}finie. On dispose donc
d'un foncteur canonique $\mbox{add}_2:\mathsf{Hmo} \rightarrow \mathsf{Hmo}_0$.
Des arguments similaires aux lemmes~\ref{addit} et ~\ref{tens} nous
permettent de montrer le lemme suivant.

\begin{lemme}
La cat{\'e}gorie $\mathsf{Hmo}_0$ est additive et le produit tensoriel
$-\otimes^{\mathbb{L}}-$ de $\mathsf{Hmo}$, voir remarque~\ref{monoi},
induit une structure mono{\"\i}dale sym{\'e}trique sur $\mathsf{Hmo}_0$.
\end{lemme} 

Soit $\mathsf{Heq}_{{kar}_0}$ la sous-cat{\'e}gorie pleine de
$\mathsf{Hmo}_0$ dont les objets sont les dg-cat{\'e}gories exactes
$\mathcal{A}$ telles que l'image du plongement
$\hat{?}:\mathrm{H}^0(\mathcal{A}) \hookrightarrow \mathcal{D}(\mathcal{A})$
est stable par facteurs directs.
\begin{remarque}
L'{\'e}quivalence entre $\mathsf{Hmo}$ et $\mathsf{Heq}_{kar}$ induit une
{\'e}quivalence entre $\mathsf{Hmo}_0$ et $\mathsf{Heq}_{{kar}_0}$.
\end{remarque}
On a donc le diagramme suivant~:
$$
\xymatrix{
\mathsf{Heq}_{{ex}_0} \ar@{^{(}->}[r]^-{\sim} & \mathsf{Hec}_0 \ar[d] \\
\mathsf{Heq}_{{kar}_0} \ar@{^{(}->}[u] \ar@{^{(}->}[r]^-{\sim} &
\mathsf{Hmo}_0\, .  
}
$$
La d{\'e}monstration du th{\'e}or{\`e}me suivant est analogue {\`a} celle du th{\'e}or{\`e}me~\ref{semiloc}. 
Soit $F:\mathsf{Hmo}\rightarrow \mathsf{C}$ un foncteur {\`a} valeurs
dans une cat{\'e}gorie additive $\mathsf{C}$.

\begin{theoreme}\label{Univ}
Les conditions suivantes sont {\'e}quivalentes~:
\begin{itemize}
\item[1)] Le foncteur $F$ est compos{\'e} d'un foncteur additif
  $\mathsf{Hmo}_0 \rightarrow \mathsf{C}$ et du foncteur canonique
  $\mathsf{Hmo}\rightarrow \mathsf{Hmo}_0$.
\item[2)] Pour toutes dg-cat{\'e}gories $\mathcal{A}$, $\mathcal{B}$, l'identit{\'e}
  $F(\left[X\right])+F(\left[Z\right])=F(\left[Y\right])$, est v{\'e}rifi{\'e}e
  dans $\mathrm{Hom}_{\mathsf{C}}(F(\mathcal{A}), F(\mathcal{B}))$ pour
  tout triangle $X \rightarrow Y \rightarrow Z \rightarrow X\left[1\right]$
  de $\mathsf{rep}_{mor}(\mathcal{A}, \mathcal{B})$.
\item[3)] Pour toute dg-cat{\'e}gorie $\mathcal{A}$, le morphisme 
$$
\xymatrix{
F(\mathcal{A})\oplus F(\mathcal{A})
\ar[rr]^-{\left[F(i_1)
        \, , \,F(i_2)\right]} && F(\mbox{T}(\mathcal{A}))
}
$$
est un isomorphisme dans $\mathsf{C}$.

\item[4)] Pout toute dg-cat{\'e}gorie pr{\'e}triangul{\'e}e $\mathcal{A}$ munie de sous-dg-cat{\'e}gories pleines pr{\'e}triangul{\'e}es
  $\mathcal{B}$ et $\mathcal{C}$ qui donnent lieu {\`a} une d{\'e}composition
  semi-orthogonale
  $\mathrm{H}^0(\mathcal{A})=(\mathrm{H}^0(\mathcal{B}),
  \mathrm{H}^0(\mathcal{C}))$, voir
  \cite{Bondal}, le morphisme
$$
F(\mathcal{B})\oplus F(\mathcal{C}) \rightarrow F(\mathcal{A})
$$
induit par les inclusions est un isomorphisme dans $\mathsf{C}$.
\end{itemize}
\end{theoreme}

\begin{remarque}
Toute {\'e}quivalence d{\'e}riv{\'e}e, voir \cite{Rickard}, \cite{Neeman}, \cite{Tilting}, donne un isomorphisme dans $\mathsf{Hmo}$ et donc dans
$\mathsf{Hmo}_0$. Cependant, il existe d'autres isomorphismes dans
$\mathsf{Hmo}_0$~: si $k$ est un corps alg{\'e}briquement clos et $A$ un
$k$-alg{\`e}bre (ordinaire) de dimension finie sur $k$ et de dimension
globale finie, alors $A$ est isomorphe {\`a} son plus grand quotient
semisimple $A/\mathrm{rad}(A)$ dans $\mathsf{Hmo}_0$
(voir \cite[2.5]{Inv}), mais dans $\mathsf{Hmo}$, $A$ est isomorphe {\`a}
$A/\mathrm{rad}(A)$ seulement si $A$ est semisimple.
\end{remarque}

\subsection*{Homologies de Hochschild, Cyclique, Negative, \ldots}
On note $\mathcal{D} \mbox{Mix}$ la cat{\'e}gorie d{\'e}riv{\'e}e de la cat{\'e}gorie
des complexes mixtes sur $k$, voir \cite{Exact}. Soit $$M:\dgcat
\longrightarrow \mathcal{D}\mbox{Mix}$$ le foncteur qui, {\`a} une petite
dg-cat{\'e}gorie $\mathcal{C}$, associe son complexe mixte
$M(\mathcal{C})$, vu comme un objet dans $\mathcal{D}\mbox{Mix}$,
voir \cite{Exact}. D'apr{\`e}s \cite{Kassel}, l'homologie cyclique de
$\mathcal{C}$ est l'homologie cyclique du complexe
mixte $M(\mathcal{C})$ et de m{\^e}me pour les autres variantes de la
th{\'e}orie (Hochschild, p{\'e}riodique, negative, \ldots). Comme dans
\cite{Inv}, on montre que le foncteur
$M$ descend {\`a} un foncteur d{\'e}fini dans la localis{\'e}e $\mathsf{Hmo}$
et que celui-ci satisfait les conditions du th{\'e}or{\`e}me~\ref{Univ}. Il se factorise donc par le foncteur $\mathcal{U}_a:\dgcat
\longrightarrow \mathsf{Hmo}_0$ et donne lieu {\`a} un foncteur additif
$$
M:\mathsf{Hmo}_0 \rightarrow \mathcal{D}\mbox{Mix}\,.
$$ 
   
\subsection*{$K$-th{\'e}orie alg{\'e}brique}
On note $\mathsf{Ho}(\mathsf{Spt})$ la cat{\'e}gorie homotopique des
spectres, voir \cite{Spectre}. Soit $$K:\dgcat \longrightarrow
\mathsf{Ho}(\mathsf{Spt})$$ le foncteur qui, {\`a} une petite dg-cat{\'e}gorie
$\mathcal{C}$, associe le spectre de $K$-th{\'e}orie de Waldhausen
\cite{Wald} associ{\'e} {\`a} la cat{\'e}gorie des dg-modules cofibrants et
parfaits dont les cofibrations sont les monomorphismes et les
{\'e}quivalences faibles les quasi-isomorphismes, voir \cite{Dugger}.
Comme il a {\'e}t{\'e} montr{\'e} dans \cite{Dugger}, les r{\'e}sultats de \cite{Wald}
impliquent que le foncteur $K$ rend inversibles les dg-foncteurs de
Morita et descend donc en un foncteur $\mathsf{Hmo} \rightarrow
\mathsf{Ho}(\mathsf{Spt})$. Le th{\'e}or{\`e}me d'additivit{\'e} de Waldhausen
\cite[1.4]{Wald} montre que ce foncteur v{\'e}rifie la condition $3)$ du
th{\'e}or{\`e}me~\ref{Univ}. Il induit donc un foncteur additif $$ K: \mathsf{Hmo}_0 \rightarrow \mathsf{Ho}(\mathsf{Spt})\, .$$

\subsection*{Vision globale}
On dispose donc du diagramme suivant~:
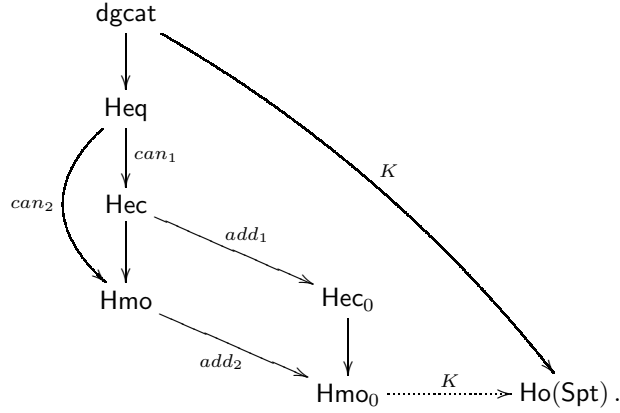
$$
\xymatrix@!0 @R=3pc @C=7pc{
 \dgcat \ar[d] 
\ar@/_-1pc/[ddddrr]^-{K} & &   \\
  \mathsf{Heq} \ar[d]^{can_1} \ar@/^-2pc/[dd]_-{can_2} & &  \\
 \mathsf{Hec} \ar[d] \ar[dr]^-{add_1} & &  \\
  \mathsf{Hmo} \ar[dr]_-{add_2} & \mathsf{Hec}_0 \ar[d] &  \\ 
  & \mathsf{Hmo}_0 \ar@{.>}[r]^-{K}  & \mathsf{Ho}(\mathsf{Spt})\,.
}
$$

\subsection*{Charact{\`e}re de Chern}
On note $\mbox{K}_0(\mathcal{C})$ le groupe de Grothendieck d'une
dg-cat{\'e}gorie $\mathcal{C}$, c'est-{\`a}-dire le groupe de Grothendieck de
la cat{\'e}gorie des objets compacts dans $\mathcal{D}(\mathcal{C})$.
Soit $\mathcal{A}$ la dg-cat{\'e}gorie avec un seul objet $3$ et telle
que $\mathrm{Hom}_{\mathcal{A}}(3,3)=k$.
\begin{lemme}
On a un isomorphisme naturel de groupes abeliens
$$\mathrm{Hom}_{\mathsf{Hmo}_0}(\mathcal{A}, \mathcal{C})
\stackrel{\sim}{\longrightarrow} \mbox{K}_0(\mathcal{C})\, .$$
\end{lemme}

\begin{proof}
En effet le groupe abelien $\mathrm{Hom}_{\mathsf{Hmo}_0}(\mathcal{A},
\mathcal{B})$ est par d{\'e}finition {\'e}gal {\`a}
$\mbox{K}_0(\mathsf{rep}_{mor}(\mathcal{A}, \mathcal{C}))$. On sait
que la cat{\'e}gorie $\mathsf{rep}_{mor}(\mathcal{A}, \mathcal{C})$
s'identifie {\`a} $\mathsf{rep}(\mathcal{A}, \mathcal{C}_{fib})$, o{\`u}
$\mathcal{C}_{fib}$ est un remplacement fibrant de la dg-cat{\'e}gorie
$\mathcal{C}$ par rapport {\`a} la structure de cat{\'e}gorie de mod{\`e}les de
Quillen d{\'e}crite dans le th{\'e}or{\`e}me~\ref{theorem2}. On remarque que la
cat{\'e}gorie $\mathsf{rep}(\mathcal{A}, \mathcal{C}_{fib})$ s'identifie {\`a}
la sous-cat{\'e}gorie pleine de $\mathcal{D}(\mathcal{C})$ form{\'e}e des
objets compacts.
\end{proof}
Soit $F:\mathsf{Hmo}_0 \rightarrow \mathrm{Mod}\,\mathbb{Z}$ un
foncteur. Le lemme de Yoneda et le lemme pr{\'e}c{\'e}dent impliquent le lemme suivant.

\begin{lemme}
On a un isomorphisme de groupes
$$ \mathrm{Nat}(\mbox{K}_0, F) \stackrel{\sim}{\longrightarrow}
F(\mathcal{A})\, ,$$
o{\`u} $\mathrm{Nat}(\mbox{K}_0, F)$ d{\'e}signe le groupe ab{\'e}lien des
transformations naturelles de $\mbox{K}_0$ vers $F$.
\end{lemme}

A titre d'exemple, pour $n \in \mathbb{N}$, soit
$\mathrm{HC}_n(\mathcal{C})$ le $n$-i{\`e}me groupe d'homologie
cyclique. On sait d{\'e}j{\`a} que le foncteur 
$$\mathrm{HC}_n: \dgcat \rightarrow \mathrm{Mod}\, \mathbb{Z}$$ se
factorise par $\dgcat \rightarrow \mathsf{Hmo}_0$. A partir de
l'isomorphisme
$$ \mathrm{HC}_*(k)\simeq k \left[ u \right], \, \, \left| u \right|
=2\, , $$
le lemme nous fournit les caract{\`e}res de Chern
$$ch_{2n}: \mbox{K}_0 \rightarrow \mathrm{HC}_{2n} \, . $$

\section{dg-cat{\'e}gories compactes et lisses}
\begin{remarque}\label{nul}
On observe facilement que la cat{\'e}gorie aditive $(\mathsf{Hec}_0,
\oplus)$ poss{\`e}de des sommes infinies. Cela implique que son groupe de
Grothendieck $\mbox{K}_0(\mathsf{Hec}_0, \oplus)$  est nul.
\end{remarque}
Supposons que $k$ est un corps. Suivant \cite{Kontsevich}, une dg-cat{\'e}gorie $\mathcal{A}$ est {\em compacte} si
\begin{itemize}
\item son ensemble d'objets est fini et
\item pour tous objets $X$ et $Y$ dans $\mathcal{A}$, la dimension de 
  $\mathrm{H}^*(\mathrm{Hom}_{\mathcal{C}}(X, Y))$ est finie;
\end{itemize}
elle est {\em lisse} si le $\mathcal{A}$-$\mathcal{A}$-bimodule associ{\'e} {\`a} l'identit{\'e}
  de $\mathcal{A}$, qu'on note $_{\mathcal{A}}\mathcal{A}_{\mathcal{A}}$, voir
  remarque~\ref{bimod}, est un objet compact dans
  $\mathcal{D}(\mathcal{A}^{op}\otimes \mathcal{A})$.
On note $\mathsf{Hec}_0^{cl}$ la sous-cat{\'e}gorie pleine de
$\mathsf{Hec}_0$ dont les objets sont les dg-cat{\'e}gories compactes et lisses.

\begin{lemme}
La cat{\'e}gorie $\mathsf{Hec}_0^{cl}$ est une sous-cat{\'e}gorie additive
mono{\"\i}dale de $\mathsf{Hec}_0$.
\end{lemme}
\begin{proof}
Pour des dg-cat{\'e}gories compactes et lisses $\mathcal{A}$ et $\mathcal{B}$, on
remarque que 
$$
\mathcal{D}((\mathcal{A}\amalg\mathcal{B})^{op}\otimes(\mathcal{A}\amalg\mathcal{B}))\backsimeq
\mathcal{D}(\mathcal{A}^{op}\otimes \mathcal{A})\times
\mathcal{D}(\mathcal{A}^{op}\otimes \mathcal{B})\times
\mathcal{D}(\mathcal{B}^{op}\otimes \mathcal{A})\times \mathcal{D}(\mathcal{B}^{op}\otimes \mathcal{B})\,.$$
Puisque le bimodule $_{\mathcal{A}\amalg\mathcal{B}}\mathcal{A}\amalg\mathcal{B}_{\mathcal{A}\amalg\mathcal{B}}$
  s'identifie {\`a} $_{\mathcal{A}}\mathcal{A}_{\mathcal{A}}\times 0 \times 0 \times
_{\mathcal{B}}\mathcal{B}_{\mathcal{B}}$, il est compact. Clairement, $\mathcal{A}\amalg\mathcal{B}$ est encore
compacte. Le fait que $\mathcal{A}\otimes\mathcal{B}$ soit lisse
r{\'e}sulte de l'isomorphisme
$$
\mathbb{R}\mathrm{Hom}_{{(\mathcal{A}\otimes\mathcal{B})}^{op}\otimes(\mathcal{A}\otimes\mathcal{B})}(\mathcal{A}\otimes\mathcal{B},
?) \backsimeq
\mathbb{R}\mathrm{Hom}_{\mathcal{A}^{op}\otimes\mathcal{A}}(\mathcal{A},
\mathbb{R}\mathrm{Hom}_{\mathcal{B}^{op}\otimes\mathcal{B}}(\mathcal{B}, ?))\,.
$$
Par le th{\'e}or{\`e}me de K{\"u}nneth, $\mathcal{A}\otimes\mathcal{B}$ est encore compacte.
\end{proof}

On note $\mathrm{HH}_0(\mathcal{A})$ le $0$-i{\`e}me groupe d'homologie
de Hochschild d'une dg-cat{\'e}gorie $\mathcal{A}$. On dispose donc d'un
foncteur additif $$\mathrm{HH}_0: \mathsf{Hec}_0 \longrightarrow
\mathrm{Mod}\, k\,.$$ En outre, si la dg-cat{\'e}gorie
$\mathcal{A}$ est compacte et lisse, alors $\mathrm{HH}_0(\mathcal{A})$ est de dimension
finie. On dispose du diagramme~:
$$
\xymatrix{
(\mathsf{Hec}_0, \oplus) \ar[r]^-{\mathrm{HH}_0} &
\mathrm{Mod}\, k \\
*+<1pc>{(\mathsf{Hec}_0^{cl}, \oplus)} \ar@{^{(}->}[u]
\ar[r]^-{\mathrm{HH}_0} & *+<1pc>{\mathrm{mod}\, k\,,}
\ar@{^{(}->}[u]
}
$$
o{\`u} $\mathrm{mod}\, k$ d{\'e}signe la cat{\'e}gorie des espaces vectoriels de dimension finie. En applicant le foncteur $\mbox{K}_0$ on
obtient 
$$
\xymatrix{
 \mathrm{K}_0(\mathsf{Hec}_0^{cl}, \oplus)
 \ar[rr]^-{\mathrm{K}_0(\mathrm{HH}_0)} & & 
\mbox{K}_0(\mathrm{mod}\,k) \,.
}
$$
En particulier, ces anneaux commutatifs, avec la structure
multiplicative induite par le produit tensoriel, sont non nuls.

Soit $\mathcal{P}\mathcal{T}$ le groupe de Grothendieck d{\'e}fini dans
\cite{Grothendieck}. Rappelons que les auteurs de \cite{Grothendieck},
\cite{Bondal} appelent {\em pr{\'e}triangul{\'e}es} les petites dg-cat{\'e}gories
$\mathcal{A}$ telles que $\mathcal{A} \hookrightarrow
\mbox{pre-tr}(\mathcal{A})$ est une quasi-{\'e}quivalence (c'est-{\`a}-dire
que $\mathcal{A}$ est fibrante pour la structure de cat{\'e}gorie de
mod{\`e}les de Quillen quasi-{\'e}quiconique d{\'e}crite dans le
th{\'e}or{\`e}me~\ref{theorem}). Le groupe $\mathcal{P}\mathcal{T}$ est par
d{\'e}finition le groupe ab{\'e}lien engendr{\'e} par les classes de
quasi-{\'e}quivalence $\left[\mathcal{A}\right]$ de petites dg-cat{\'e}gories
pr{\'e}triangul{\'e}es soumises aux relations qui proviennent des
d{\'e}compositions semi-orthogonales.
\par
L'argument de la remarque~\ref{nul}
montre qu'en fait le groupe $\mathcal{P}\mathcal{T}$ s'annule. Pour
obtenir un groupe non nul, il convient d'imposer des conditions de
finitude~: soit ${\mathcal{P}\mathcal{T}}^{cl}$ le groupe ab{\'e}lien
engendr{\'e} par les classes de quasi-{\'e}quivalences $\left[ \mathcal{A}
\right]$ de petites dg-cat{\'e}gories $\mathcal{A}$ compactes, lisses et
pr{\'e}triangul{\'e}es soumises aux relations qui proviennent des
d{\'e}compositions semi-orthogonales. Le produit tensoriel munit ${\mathcal{P}\mathcal{T}}^{cl}$ d'une structure d'anneau commutatif.

\par
Par exemple, si $X$ est une vari{\'e}t{\'e} projective lisse, alors la
cat{\'e}gorie d{\'e}riv{\'e}e $\mathcal{D}^b(\mathsf{coh}X)$ est {\'e}quivalent {\`a}
$\mathrm{H}^0(\mathcal{D}^b(\mathsf{coh}X)_{dg})$ pour une
dg-cat{\'e}gorie compacte, lisse, pr{\'e}triangul{\'e}e
$\mathcal{D}^b(\mathsf{coh}X)_{dg}$ canonique {\`a} isomorphisme pr{\`e}s dans
$\mathsf{Hmo}$. Par exemple, pour $\mathcal{D}^b(\mathsf{coh}X)_{dg}$,
on peut prendre la dg-cat{\'e}gorie des complexes born{\'e}s {\`a} gauche de
$\mathcal{O}_X$-modules injectifs dont l'homologie est born{\'e}e et
coh{\'e}rente.

\begin{proposition}\label{surj}
On a un morphisme surjectif d'anneaux commutatifs
$$
{\mathcal{P}\mathcal{T}}^{cl} \rightarrow \mbox{K}_0(\mathsf{Hec}_0^{cl}, \oplus)\, .
$$
\end{proposition}

\begin{proof}

On note $\left[\mathcal{A}\right]$ la classe de quasi-{\'e}quivalence
d'une dg-cat{\'e}gorie compacte, lisse et pr{\'e}triangul{\'e}e et $\left[\left[
    \mathcal{A} \right]\right]$ la classe d'isomorphisme de
$\mathcal{A}$ dans $\mathsf{Hec}_0^{cl}$. On consid{\`e}re l'application qui, {\`a} $\left[\mathcal{A}\right]$, associe
$\left[\left[ \mathcal{A} \right]\right]$. Elle est clairement bien
d{\'e}finie. Pour qu'elle induise un morphisme de
${\mathcal{P}\mathcal{T}}^{cl}$ vers $\mbox{K}_0(\mathsf{Hec}_0^{cl}, \oplus)$, il reste
{\`a} verifier l'{\'e}galit{\'e}
$$
\left[\left[ \mathcal{B}\oplus \mathcal{C} \right]\right]=\left[\left[\mathcal{A}\right]\right]\, ,
$$
pour toutes dg-cat{\'e}gories $\mathcal{A}$, $\mathcal{B}$ et
$\mathcal{C}$ qui satisfont les relations dans la d{\'e}finition de
${\mathcal{P}\mathcal{T}}^{cl}$. On consid{\`e}re le diagramme
$$
\xymatrix{
 &&&& *+<1pc>{\mathcal{A}} \ar@{^{(}->}[d]^-{inc} \\
\mathcal{B}\oplus\mathcal{C}=\mbox{pre-tr}(\mathcal{B}\coprod \mathcal{C})
\ar[rrrr]^-{pre-tr(\left[ inc(\mathcal{B}) \, \,
    inc(\mathcal{C})\right])} &&&& \mbox{pre-tr}(\mathcal{A})\,.
}
$$
Le dg-foncteur $inc$ est une quasi-{\'e}quivalence et il reste {\`a} montrer que l'image du dg-foncteur
$\mbox{pre-tr}(\left[ inc(\mathcal{B}) \,\, inc(\mathcal{C})\right])$
dans $\mathsf{Hec}_0^{cl}$ est un isomorphisme. Par le lemme de
Yoneda, il s'agit de montrer que, pour toute dg-cat{\'e}gorie
$\mathcal{U}$, l'application
$$
\mbox{K}_0(\mathsf{rep}_{ec}(\mathcal{U},
\mbox{pre-tr}(\mathcal{B}\amalg \mathcal{C}))) \longrightarrow \mbox{K}_0(\mathsf{rep}_{ec}(\mathcal{U},
\mbox{pre-tr}(\mathcal{A})))
$$
est bijective. En effet, la d{\'e}composition semi-orthogonale
$\mbox{H}^0(\mathcal{A})=(\mbox{H}^0(\mathcal{B}),
\mbox{H}^0(\mathcal{C}))$, induit une
d{\'e}composition semi-orthogonale
$$
\mathsf{rep}_{ec}(\mathcal{U},
\mbox{pre-tr}(\mathcal{A}))=(\mathsf{rep}_{ec}(\mathcal{U},
\mbox{pre-tr}(\mathcal{B})), \mathsf{rep}_{ec}(\mathcal{U},
\mbox{pre-tr}(\mathcal{C}))) \, .
$$
La surjectivit{\'e} du morphisme
${\mathcal{P}\mathcal{T}}^{cl}\rightarrow
\mbox{K}_0(\mathsf{Hec}_0^{cl}, \oplus)$ est
claire par construction.
\end{proof}

\begin{remarque}
En particulier, la proposition~\ref{surj} montre que
${\mathcal{P}\mathcal{T}}^{cl}$ est non nul. Si on d{\'e}finit
${\mathcal{P}\mathcal{T}}^{cl}_{kar}$ de fa{\c c}on analogue {\`a} partir des
dg-cat{\'e}gories compactes, lisses, pr{\'e}triangul{\'e}es $\mathcal{A}$ telles
que $\mathrm{H}^0(\mathcal{A})$ est Karoubienne, on obtient la
surjection
$$
{\mathcal{P}\mathcal{T}}^{cl}_{kar} \rightarrow \mbox{K}_0(\mathsf{Hmo}_0^{cl})
$$
mentionn{\'e}e dans l'introduction.
\end{remarque}

\section*{Remerciements}
Ce travail fait partie de ma th{\'e}se sous la
direction de B. Keller, que je remercie pour de nombreuses conversations
utiles.

\label{}



\begin{thebibliography}{00}

%


\bibitem{Tilting} L.~Angeleri-H{\"u}gel, D.~Happel, H.~Krause (eds.),
  Handbook of tilting theory, to appear.

\bibitem{Balmer} P.~Balmer, M.~Schlichting, {\em Idempotent completion
    of triangulated categories}, J.~Alg. {\bf 236} (2001), no.~2, 819--834. 

\bibitem{Bondal} A.~Bondal, M.~Kapranov, {\em Framed triangulated
    categories} (Russian) Mat.~Sb. {\bf 181} (1990) no.~5, 669--683;
    translation in Math.~USSR-Sb. {\bf 70} no.~1, 93--107.

\bibitem{Grothendieck} A.~Bondal, M.~Larsen, V.~Lunts, {\em
    Grothendieck ring of pretriangulated categories},
    Int.~Math.~Res.~Not.~(2004), no.~29, 1461--1495.

\bibitem{Bon-Orl} A.~Bondal, D.~Orlov, {\em Semiorthogonal
    decomposition for algebraic varieties}, preprint MPIM 95/15
    (1995), preprint math. AG/9506012.

\bibitem{Drinfeld} V.~Drinfeld, {\em DG quotients of DG categories},
J. Alg. {\bf 272} (2004), 643--691.

\bibitem{Drinfeldc} V.~Drinfeld, {\em DG categories}, Series of talks
  at the Geometric Langlands Seminar, Chicago, Notes taken by
  D.~Ben-Zvi, Fall 2002.

\bibitem{Dugger} D.~Dugger, B.~Shipley, {\em K-theory and derived
    equivalences}, Duke Math.~J. {\bf 124} (2004), no.~3, 587--617.

\bibitem{Spectre} M.~Hovey, B.~Shipley, J.~Smith, {\em Symmetric
  spectra}, J.~Amer.~Math.~Soc. {\bf 13} (2000), no.~1, 149--208.

\bibitem{Hovey} Mark Hovey, {\em Model Categories}, Mathematical
  Surveys and Monographs, volume 63, AMS, Providence, RI, 1999.

\bibitem{Kassel} C.~Kassel, {\em Cyclic homology, comodules and mixed
    complexes}, J.~Algebra {\bf 107} (1987), no.~1, 195--216.

\bibitem{Orbit} B.~Keller, {\em On triangulated orbit categories},
preprint, math.~RT/0503240.

\bibitem{Exact} B.~Keller, {\em On the cyclic homology of exact
    categories}, J.~Pure~Appl.~Algebra {\bf 136} (1999), no.~1, 1--56.

\bibitem{Inv} B.~Keller, {\em Invariance and Localization for cyclic
    homology of DG algebras}, J.~Pure~Appl.~Algebra {\bf 123} (1998),
    no.~1--3, 223--273.

\bibitem{Keller} B.~Keller, {\em Deriving DG categories},
  Ann. Scient. Ec. Norm. Sup. {\bf 27} (1994), 63--102.

\bibitem{Kontsevich} M.~Kontsevich, {\em Topological field theory for
    triangulated categories}, Talk at the conference on $K$-theory and
    Noncommutative Geometry, Institut Henri Poincar{\'e}, Paris, June 2004. 

\bibitem{Kontsevichc} M.~Kontsevich, {\em Triangulated categories and
    geometry}, Course at the {\'E}cole Normale Sup{\'e}rieure, Paris, Notes
    taken by J.~Bella{\"\i}che, J.-F.~Dat, I.~Marin, G.~Racinet and
    H.~Randriambololona, 1998.

\bibitem{Neeman} A.~Neeman, {\em The K-theory of triangulated
    categories}, in: Handbook of $K$-theory, edited by
    E.~Friedlander, D.~Grayson, Springer, Berlin, 2005.

\bibitem{Quillen} D.~Quillen, {\em Homotopical algebra}, Lecture Notes
  in Mathematics, {\bf 43}, Springer-Verlag, 1967.

\bibitem{Rickard1} J.~Rickard, {\em Derived equivalences as derived
    functors}, J.~London Math.~Soc.~(2) {\bf 43} (1991), no.~1, 37--48.

\bibitem{Rickard} J.~Rickard, {\em Morita theory for derived
  categories}, J.~London Math.~Soc. (2), {\bf 39} (1989), no.~3, 436--456.

\bibitem{Rouquier} R.~Rouquier, A.~Zimmermann, {\em Picard groups for
    derived module categories}, Proc.~London Math.~Soc.~(3) {\bf 87}
    (203), no.~1, 197--225.

\bibitem{cras} G.~Tabuada, {\em Une structure de cat{\'e}gorie de mod{\`e}les
    de Quillen sur la cat{\'e}gorie des dg-cat{\'e}gories},
    C. R. Math. Acad. Sci. Paris {\bf 340} (2005) no.~1, 15--19. 

\bibitem{Toen} B.~To{\"e}n, {\em The homotopy theory of dg-categories and
    derived Morita theory}, preprint, math.~AG/0408337.

\bibitem{Wald} F.~Waldhausen, {\em Algebraic K-theory of spaces},
  Algebraic and geometric topology (New Brunswick, N.~J., 1983),
  318--419, Lecture Notes in Math., 1126, Springer, Berlin, 1985.

\end{thebibliography}
\end{document}